\documentclass[11pt]{article}
\usepackage{}
\usepackage{amsfonts}
\usepackage{amssymb}
\usepackage{amsmath}

\def\sqr#1#2{{\vcenter{\vbox{\hrule height.#2pt
              \hbox{\vrule width.#2pt height#1pt \kern#1pt \vrule
width.#2pt}
              \hrule height.#2pt}}}}
\def\signed #1{{\unskip\nobreak\hfil\penalty50
              \hskip2em\hbox{}\nobreak\hfil#1
              \parfillskip=0pt \finalhyphendemerits=0 \par}}
\def\endpf{\signed {$\sqr69$}}
\def\dbN{{\mathbb{N}}}

\def\3n{\negthinspace \negthinspace \negthinspace }
\def\2n{\negthinspace \negthinspace }
\def\1n{\negthinspace }

\def\ds{\displaystyle}
\def\D{\Delta}

\def\no{\noindent}

\def\ms{\medskip}
\def\bs{\bigskip}
\def\q{\quad}
\def\qq{\qquad}

\def\div{\mathop{\nabla\cd\neg}}

\def\cd{\cdot}

\def\deq{\mathop{\buildrel\D\over=}}

\def\({\Big (}
\def\){\Big )}
\def\[{\Big[}
\def\]{\Big]}
\def\neg{\negthinspace}

\def\be{\begin{equation}}
\def\bel{\begin{equation}\label}
\def\ee{\end{equation}}
\def\bea{\begin{eqnarray}}
\def\eea{\end{eqnarray}}
\def\bt{\begin{theorem}}
\def\et{\end{theorem}}
\def\bc{\begin{corollary}}
\def\ec{\end{corollary}}
\def\bl{\begin{lemma}}
\def\el{\end{lemma}}
\def\bp{\begin{proposition}}
\def\ep{\end{proposition}}
\def\br{\begin{remark}}
\def\er{\end{remark}}
\def\ba{\begin{array}}
\def\ea{\end{array}}
\def\bd{\begin{definition}}
\def\ed{\end{definition}}

\newtheorem{lemma}{Lemma}[section]
\newtheorem{remark}{Remark}[section]

\newtheorem{theorem}{Theorem}[section]
\newtheorem{corollary}{Corollary}[section]

\newtheorem{definition}{Definition}[section]
\newtheorem{proposition}{Proposition}[section]

\makeatletter
   
   \@addtoreset{equation}{section}
\makeatother

\begin{document}

\title{\bf   A  Class of Multi-objective Control problems for   Quasi-linear  Parabolic Equations  \thanks{This work was
 partially supported
by National Key R\&D Program of China under grant 2023YFA1009002,  
NSF of China under grants 12371444 and 11931011, New Cornerstone Investigator Program, and the Science Development Project of Sichuan University under grant 2020SCUNL201.}}

\author{Yanming Dong,\thanks{School of Mathematics and Statistics, Northeast Normal
University, Changchun 130024, China.  E-mail address:
dongym946@nenu.edu.cn.}\quad     Xu Liu\thanks{School of Mathematics and Statistics, Northeast Normal
University, Changchun 130024, China.   E-mail: 
liux216@nenu.edu.cn.}   \ and\  Xu Zhang\thanks{School of Mathematics,  Sichuan University, 
Chengdu  610064, China. E-mail address: zhang\_xu@scu.edu.cn.}}

\date{}

\maketitle

\begin{abstract} 
This paper is devoted to studying a multi-objective control problem  for a class of  
multi-dimensional  quasi-linear  parabolic equations.
The considered   system is driven by  a  leader  control  and two follower controls. 
For each leader  control,  a pair of follower controls is searched for as a Nash quasi-equilibrium
(or  Nash equilibrium)  of  cost functionals,    while 
the aim for a leader  control is to solve a  controllability problem. 
This hierarchic control  problem  may be transformed into
 controllability of 
a strongly coupled system of quasi-linear parabolic equations through  one control.  
 Regarding   controllability   for    quasi-linear parabolic equations of second order,  
  the existing   results usually  require  coefficients in principal parts  to be independent of gradient of solutions,
   or spacial dimension to be limited.   In this paper,  the coefficients in 
 principal parts for the controlled  quasi-linear system contain  not  only  the state itself
   but also   its gradient with general 
 spacial dimension.   
 \end{abstract}

\bs

\no{\bf Key Words}.   Quasi-linear parabolic equation,  controllability,
  Stackelberg-Nash  strategy, Carleman estimate

\ms

\no{\bf AMS subject classifications}. 93B05,  93C20,  35K59

\date{}
\maketitle

\section{Introduction and main results }

Let  $T>0$ and  $\Omega $ be a   bounded domain in 
$\mathbb{R}^n$ $(n\in\mathbb N)$ with a smooth boundary $\partial\Omega$. 
For $j=0, 1, 2$,  assume that $\omega_j$ and $\tilde\omega_j$  are 
nonempty open subsets of $\Omega$ such that $\overline{\tilde\omega_j}\subseteq \omega_j$, 
where 
$\overline{\tilde\omega_j}$ denotes the closure of the set $\tilde\omega_j$.  Denote by $\xi_{j}\in C_0^\infty(\omega_j)$  
 smooth functions satisfying that  $\xi_{j}(x)=1$
 in $\tilde\omega_j$ and $0\leq \xi_j(x)\leq 1$ in $\omega_j$.  Put $Q=\Omega\times(0,T)$ and $\Sigma=\partial\Omega\times(0,T)$.  
 Consider the following controlled  quasi-linear 
parabolic equation:
\begin{eqnarray}\label{a}
\left\{
\begin{array}{ll}
y_t-\sum\limits^n_{i,j=1}\big(a^{ij}(y,\nabla y)y_{x_i}\big)_{x_j}+f(y, \nabla y)
=\xi_{0} u+\xi_{1} v_1+\xi_{2} v_2     &\mbox{ in }   Q,\\[2mm]
y=0  &\mbox{ on }  \Sigma,\\[2mm]
y(x,  0)=y_0(x)  &\mbox{ in }  \Omega,
\end{array}
\right.
\end{eqnarray}
where  $(u,  v_1, v_2)$ is a  triple  of control variables,  $y$ is the state variable, 
$y_0$ is an initial  value,
$f:  \mathbb R\times\mathbb R^n\rightarrow \mathbb R$ is a 
$C^4$ function with  $f(0,  {\bf  0})=0$
and $a^{ij}:\mathbb{R}\times\mathbb{R}^n\rightarrow\mathbb{R}$ are    $C^4$
functions with $a^{ij}=a^{ji}$  in  $\mathbb R^{n+1}$ for $i,j=1,\ldots,n$. Moreover,   for a  positive constant $\rho_0$, 
\begin{equation*}\label{c}
\sum_{i, j=1}^n  a^{i j}(s, \eta)\zeta_i\zeta_j\geq\rho_0|\zeta|^2,
\quad\forall\ (s, \eta, \zeta)=(s, \eta^1, \cdots, \eta^n, \zeta_1, \cdots, \zeta_n)\in\mathbb R^{1+2n}.
\end{equation*}

Quasi-linear parabolic equations  may  describe
 a wide range of diffusion phenomena in nature, 
 such as filtration, dynamics of biological groups and heat conduction. 
  For example,  suppose
  that a  
  body occupies the domain  $\Omega$ and it is anisotropic.
  If $y$, $c$ and $\rho$  denote  temperature distribution, specific heat and 
  density  
   of the 
  body, respectively, 
 the conservation law  will  give us
$$
c\rho y_t+\text{div} {\bf J}=f,
$$
where  ${\bf J}$ denotes a heat flux vector and $f$ is a heat source. 
By the Fourier experiment law,  the heat  flux vector has the form of   ${\bf J}=
-\mathcal{A}\nabla y$ 
for thermal conductivity matrix $\mathcal{A}=\big(a^{i j}\big)_{1\leq i, j\leq n}$.

When $\mathcal{A}$ and  $f$ depend  on 
 both  temperature distribution $y$ and its  rate of change (e.g. Perona-Malik 
 model and p-Laplacian equation),    
   it follows that 
$$
c \rho y_t-\sum\limits^n_{i,j=1}\big(a^{ij}(y,\nabla y)y_{x_i}\big)_{x_j}=f(y, \nabla y).
$$
For simplicity, assume $c\rho\equiv1$.
If one can influence  heat diffusion through  heat  sources $u$,  $v_1$ and $v_2$ 
at different subdomains, 
 the associated  system 
may be expressed as the form of (\ref{a}). 
The functions $\xi_j$ $(j=0, 1, 2)$ placed before  the controls $u, v_1$ and $v_2$
indicate  that they are active in 
different local domains $\omega_j$. Smoothness 
of $\xi_j$
  ensures that arguments for the quasi-linear parabolic  system
  may be conducted in the framework of classical solutions.


Next,  we  introduce some   function spaces, which will be used later.  For any $k,
\ell\in \dbN$,  denote by $C^{k,\ell}(\overline{Q})$ the set of
 functions defined in $\overline{Q}$, which have continuous derivatives up to order $k$ with
respect to the spacial  variable and up to order $\ell$ with respect to
the time variable,  and by $C^k(\overline{\Omega})$ the set of 
functions  defined in $\overline{\Omega}$,  which have continuous derivatives 
up to order $k$.  For  $\alpha\in (0, 1)$,  set
 $$
 C^{k+\alpha}(\overline{\Omega})=\Big\{\   v\in C^k(\overline{\Omega})\   \Big|\   
 \sup\limits_{x_1, x_2\in\overline{\Omega}, x_1\neq x_2}\sup\limits_{|\sigma|=k}
 \displaystyle\frac{|\partial^{\sigma}_{x}
 v(x_1)-\partial^\sigma_{x}
v(x_2)|}{|x_1-x_2|^\alpha}<\infty\    \Big\},
$$
and denote by $C_0^{k+\alpha}(\Omega)$  the set of all functions  in $C^{k+\alpha}(\overline{\Omega})$,  which have compact support in $\Omega$. Here and hereafter,   $\sigma=(\sigma_1, \cdots, \sigma_n)$ denotes a multi-index with the order $|\sigma|
=\sigma_1+\cdots+\sigma_n$.
Denote
$$\displaystyle C^{k+\alpha, \frac{k+\alpha}{2}}(\overline{Q})=\Big\{\  v
\in C(\overline{Q})\  \Big|\   \partial_x^\sigma\partial_t^r v\in C(\overline{Q})
\  \mbox{ and }\ [v]_{k, \alpha}<\infty, 
\mbox{ for }  0\leq 2r+|\sigma|\leq k\  \Big\},
$$
 and  for an even number $k$,
$$
\displaystyle [v]_{k, \alpha}=\sum\limits_{2r+|\sigma|=k} 
\sup\limits_{(x_1,t_1)\neq(x_2,t_2)}
\frac
{|\partial^\sigma_x\partial^{r}_t
v(x_1,t_1)-\partial^\sigma_x\partial_t^{r}
v(x_2,t_2)|}{(|x_1-x_2|+|t_1-t_2|^{1/2})^\alpha}
$$
and for an odd number $k$,
\begin{eqnarray*}
&&\displaystyle [v]_{k, \alpha}=\sum\limits_{2r+|\sigma|=k} 
\sup\limits_{(x_1,t_1)\neq(x_2,t_2)}
\frac
{|\partial^\sigma_x\partial^{r}_t
v(x_1,t_1)-\partial^\sigma_x\partial_t^{r}
v(x_2,t_2)|}{(|x_1-x_2|+|t_1-t_2|^{1/2})^\alpha}\\
&& \quad\quad\quad\quad+
\sum\limits_{2r+|\sigma|=k-1} \sup\limits_{x\in\overline{\Omega}}
\sup\limits_{t_1\neq t_2}
\frac
{|\partial^\sigma_x\partial^{r}_t
v(x,t_1)-\partial^\sigma_x\partial_t^{r}
v(x,t_2)|}{|t_1-t_2|^{\frac{1+\alpha}{2}}},
\end{eqnarray*}
with the norm 
$$|v|_{C^{k+\alpha, \frac{k+\alpha}{2}}(\overline{Q})}=
\sum\limits_{0\leq 2r+|\sigma|\leq k} \sup\limits_{(x, t)\in\overline{Q}}|\partial_x^\sigma\partial_t^r v(x, t)|+
 [v]_{k, \alpha}.
$$
 The definitions of the above function spaces can be found in \cite{La}.
Moreover,   define the space
\begin{eqnarray*}
&&\displaystyle \widetilde C^{\alpha, \frac{\alpha}{2}}_1 (\overline{Q})=\Big\{\  v
\in C^{\alpha, \frac{\alpha}{2}}(\overline{Q})\  \Big|\  
\nabla_x v\in \big(C^{\alpha, \frac{\alpha}{2}}(\overline{Q})\big)^n\ \Big\}
\end{eqnarray*}
with the norm
$$
|v|_{\widetilde C^{\alpha, \frac{\alpha}{2}}_1 (\overline{Q})}
=|v|_{C^{\alpha, \frac{\alpha}{2}}(\overline{Q})}
+\sum_{i=1}^n 
|v_{x_i}|_{C^{\alpha, \frac{\alpha}{2}}(\overline{Q})},
$$ 
where $\nabla_x v$ 
denotes the gradient of $v=v(x, t)$ with respect to the spacial  variable $x=(x_1, \cdots, x_n)$.
When it does not cause confusion, we use $\nabla v$ to denote
$\nabla_x v$.

By  well-posedness results of quasi-linear parabolic equations, there exists a positive constant 
$\rho_1$,  such that for any $y_0\in C^{2+\alpha}_0(\Omega)$  and 
$u, v_1, v_2\in C^{\alpha, \frac{\alpha}{2}}(\overline{Q})$   with 
$$
y_0\in \tilde B_{\rho_1}\deq
\Big\{\  v\in C_0^{2+\alpha}(\Omega)\   \Big|\ 
|v|_{C^{2+\alpha}(\overline{\Omega})}<\rho_1\  \Big\}$$
 and 
 $$
 u,  v_1,  v_2\in B_{\rho_1}\deq\Big\{\  v\in C^{\alpha, \frac{\alpha}{2}}(\overline{Q})\   \Big|\ 
 |v|_{C^{\alpha, \frac{\alpha}{2}}(\overline{Q})}<\rho_1\  \Big\},
$$
the equation (\ref{a}) admits a unique solution $y\in C^{2+\alpha, 1+\frac{\alpha}{2}}(\overline{Q})$. 
Moreover, 
\begin{equation}\label{XX1}
|y|_{C^{2+\alpha, 1+\frac{\alpha}{2}}(\overline{Q})}\leq C\Big(|y_0|_{C^{2+\alpha}(\overline{\Omega})}
+|u|_{C^{\alpha, \frac{\alpha}{2}}(\overline{Q})}
+|v_1|_{C^{\alpha, \frac{\alpha}{2}}(\overline{Q})}
+|v_2|_{C^{\alpha, \frac{\alpha}{2}}(\overline{Q})}\Big).
\end{equation}
Here and hereafter, we use $C$ to denote a generic positive  
constant, depending only on $n$, $T$, $\Omega$, $\alpha$, $a^{i j} $ and $f$, 
which may be different in different places.

For each leader control $u\in B_{\rho_1}$ and $k=1, 2$,   define the following 
quadratic cost  functionals:
\begin{equation}\label{d}
J_k(v_1,  v_2; u)=\frac{\mu_k}{2}  \int^T_0\int_{\omega_k}|v_k(x, t)|^2dxdt+
\frac{\nu_k}{2}  \int_Q \xi_* |y(x, t)-y_{k, d}(x, t)|^2dxdt,  \quad     
\  \forall\  v_1, v_2\in B_{\rho_1},
\end{equation}
where  $\mu_k$ and $\nu_k$  are given positive constants,   $ \mathcal{\omega}$ and $\omega'$ are two
 open subsets of $\Omega$ with $\overline{\omega'}\subseteq \omega$,  
$\xi_{*}\in C^\infty_0(\omega)$, $\xi_*=1$ in $\omega'$, 
$y_{k,d}\in  C^{1+\alpha, \frac{1+\alpha}{2}}(\overline{Q})$ are  given functions,  and 
$y$  is the solution to (\ref{a}) corresponding  to   $(u,   v_1,    v_2)\in  (B_{\rho_1})^3$. 
Note that if one hopes the state $y$ to be close to the target $y_{k ,d}$ on the whole domain 
$\Omega$,  we  take $\xi_*\equiv1$ on $\Omega$. 
In such a case,    assume that $y_{k, d}(x, T)=
0$ on $\partial\Omega$. 
Because 
the hierarchic control
 problem 
 is studied  in the framework of  classical solutions,  the above 
assumptions  may guarantee  compatibility conditions for the systems involved in arguments  to hold.

For a leader control $u\in  B_{\rho_1}$,  
the action of the associated follower controls $(v_1, v_2)$ is to reach  a Nash equilibrium or Nash  quasi-equilibrium 
of the functionals $J_1$  and $J_2$.
We first recall the notions of  Nash equilibrium and Nash  quasi-equilibrium.
\begin{definition}\label{+1}
For a   leader control   $u\in B_{\rho_1}$,   
a  pair $(\bar{v}_1, \bar{v}_2)\in  (B_{\rho_1})^2$  is called a Nash equilibrium of 
the functionals $J_1$  and  $J_2$,  if 
\begin{equation}\label{+2}
  J_1(\bar{v}_1, \bar{v}_2; u)\leq  J_1(v_1,  \bar{v}_2; u),    \quad \forall\   v_1\in 
  B_{\rho_1},
\end{equation}
and
\begin{equation}\label{+3}
  J_2(\bar{v}_1,  \bar{v}_2;  u)\leq  J_2(\bar{v}_1, v_2; u),    \quad \forall\  v_2\in 
  B_{\rho_1}.
\end{equation}

On the other hand, a pair $(\bar{v}_1,\bar{v}_2)\in  (B_{\rho_1})^2$  is called a 
Nash quasi-equilibrium of $J_1$  and  $J_2$, 
 if 
\begin{equation}\label{+4}
J_{1, v_1}(\bar{v}_1, \bar{v}_2;  u)v_1
  \deq\lim\limits_{\epsilon\rightarrow 0}\frac{J_1(\bar v_1+\epsilon v_1, \bar v_2; u)
  -J_1(\bar v_1, \bar v_2; u)}{\epsilon}=0,
   \quad \forall\   v_1\in 
  B_{\rho_1},
\end{equation}
and
\begin{equation}\label{+5}
J_{2, v_2}(\bar{v}_1, \bar{v}_2;  u)v_2
 \deq\lim\limits_{\epsilon\rightarrow 0}\frac{J_2(\bar v_1, \bar v_2+\epsilon v_2; u)-
 J_2(\bar v_1, \bar v_2; u)}{\epsilon}=0, 
 \quad \forall\   v_2\in 
  B_{\rho_1}.
\end{equation}
\end{definition}

Nash equilibrium refers to a combination of strategies  in a game.
Any player who changes the strategy  will not increase their own profits, 
when other players' strategies remain unchanged.
In the previous game problem with the cost functionals $J_1$ and $J_2$, 
the profits  are to pay less control cost and to make  the associated state fully close to the desired goal. 
Note that a Nash equilibrium of $J_1$ and $J_2$ is also a  Nash quasi-equilibrium,
while the opposite may not necessarily be true. 
But when  both $J_1(\cdot, \bar v_2; u)$  and  $J_2(\bar v_1, \cdot; u)$ are  
convex  in $B_{\rho_1}$,   
 the Nash quasi-equilibrium  $(\bar{v}_1, \bar{v}_2)$ becomes a Nash equilibrium. 
 In order to guarantee $(\bar{v}_1, \bar{v}_2)$ to be a Nash equilibrium, 
 more requirements in general are  needed for the system (\ref{a}) and the cost functionals 
 (see Section \ref{subsection}). 
 
 This paper is devoted to studying a multi-objective control problem 
 for the quasi-linear parabolic equation (\ref{a}).
  For any 
 leader control $u\in  B_{\rho_1}$, which is active on the local domain  
 $\omega_0$,   it  will be proven  that under suitable conditions,   the functionals $J_1(\cdot, \cdot;  u)$  
 and $J_2(\cdot, \cdot;  u)$
 have a unique  Nash equilibrium or Nash quasi-equilibrium.   Next,   
 a leader control $\bar u\in  B_{\rho_1}$  will be 
 found so that the corresponding solution  $y(\cdot, \cdot;  
 \bar u,  \bar  v_1,   \bar v_2)$ to (\ref{a}) reaches the desired target zero at   given time $T$,  i. e., 
 $$
 y(\cdot, T;  \bar u,  \bar  v_1,   \bar v_2)=0\quad\text{  in  }\Omega,
 $$
 where $(\bar v_1, \bar v_2)\in(B_{\rho_1})^2$ is the Nash equilibrium or 
 Nash quasi-equilibrium for the leader control $\bar u$.

  In  the above control problem,   
 the triple of controls are divided into  two levels.  
 The leader control is dominant and
  the decision of the  follower controls 
 will be influenced by the leader one.    
 For example,     strategy  of   government in a city may be viewed as a leader control and 
different decisions of some real estate agencies  serve as follower controls. 
When the government issues a new policy that affects housing prices, 
the real estate agencies  engage in a game.  
This multi-objective control problem is to find a strategy  of  the government,  according to 
which,   the decisions of these agencies  can reach a Nash equilibrium in the game. Moreover,
 the behaviors  of the government and agencies 
 can enable the housing price to reach an expected target.

 Up to now,  there are numerous works on  controllability  
   of partial differential equations by hierarchic 
 controls.   We refer to \cite{Al,  Ar, Ara, Ca, Car, Djom, Ga, H, Lion, Lio} and  references therein for some known results  in this respect.  
 For  quasi-linear parabolic equations,
 in \cite{N} and \cite{Ro}  the  coefficients in principal parts of the considered   equations
 did not contain  gradient of  solutions 
 and  spatial dimension was limited to no more than three.
 The limitation on spatial dimension was removed in \cite{Ni}.
 In \cite{Nu},  the  coefficients in principal parts of the quasi-linear parabolic equations
 depended on  gradient of  solutions, but 
 only the problems  in 
  one  spacial dimension  were  analyzed.
 This paper is devoted to 
  studying the multi-objective control problem for 
  the 
  quasi-linear parabolic equation (\ref{a}), in which, 
 the   coefficients of    principal part  depend on 
    both  solution and its gradient, and  spacial dimension is general.

  In order to get the desired controllability  result  for 
  (\ref{a})   through  hierarchic  controls,  we impose 
 the following assumptions: 
  \begin{eqnarray}\label{assume1}
 \begin{array}{ll}
 \displaystyle (1)\  \tilde\omega_0\cap \omega' \neq \emptyset;  \mbox{and }&\\[2mm]
 \displaystyle (2)\   \mbox{in }(\ref{d}),   
\  \mbox{the constants }\mu_k\geq \mu_0\mbox{ with }\mu_0\mbox{ being the constant in } (\ref{mu**})
 \mbox{ for }k=1, 2. &
  \end{array}
 \end{eqnarray}
 
 \medskip

The first  main result  of this paper is stated as follows.
\begin{theorem}\label{th1}
Assume that the conditions in $(\ref{assume1})$ hold.   There exists a positive 
constant $\delta_0$, such that for any  initial value $y_0\in C_0^{3+\alpha}(\Omega)$ and $y_{k, d}\in 
 C^{1+\alpha, \frac{1+\alpha}{2}}(\overline{Q})$
satisfying that 
$$
|y_0|_{C^{3+\alpha}(\overline{\Omega})}
+\sum_{k=1}^2|y_{k, d}|_{C^{1+\alpha, \frac{1+\alpha}{2}}(\overline{Q})}
+\sum_{k=1}^2|\hat\rho y_{k, d}|_{L^2(Q)}\leq \delta_0,
$$ 
one  can always  find  a leader  control
 $\bar u\in C^{1+\alpha, \frac{1+\alpha}{2}}(\overline{Q})$ so that    $J_1(\cdot, \cdot; \bar u)$  
and  $J_2(\cdot, \cdot; \bar u)$ have a unique  pair of 
 Nash quasi-equilibrium $(\bar{v}_1,\bar{v}_2)\in (B_{\rho_1})^2$.  Moreover,   the corresponding  solution 
 $y(\cdot, \cdot; \bar u, \bar{v}_1,\bar{v}_2)$ to $(\ref{a})$
satisfies that $$y(\cdot,T; \bar u, \bar{v}_1,\bar{v}_2)=0\quad \mbox{   in } \Omega,$$
where $\hat\rho$ is a function given in $(\ref{notation})$.
\end{theorem}

Furthermore,  for any leader control, 
 in order to derive the existence of  a pair of Nash equilibrium,  
 we need more condition and assume that
  \begin{eqnarray}\label{assume2}
 \begin{array}{ll}
 \displaystyle    (3) \   \mbox{in }(\ref{d}), 
\  \mu_k\geq \mu^0\mbox{ with }\mu^0\mbox{ being the constant in } (\ref{mu****})
 \mbox{ for }k=1, 2.&
 \end{array}
 \end{eqnarray}
In this case,  the associated   controllability result for (\ref{a}) by hierarchic  controls 
 is given as follows. 
\begin{theorem}\label{+6}
Assume that the conditions in $(\ref{assume1})$ and $(\ref{assume2})$ hold.   
There exists a positive constant $\delta^0$, such that for any  initial value $y_0\in C_0^{3+\alpha}(\Omega)$ and $y_{k, d}\in C^{1+\alpha, \frac{1+\alpha}{2}}(\overline{Q})$
satisfying that 
$$
|y_0|_{C^{3+\alpha}(\overline{\Omega})}
+\sum_{k=1}^2 |y_{k, d}|_{C^{1+\alpha, \frac{1+\alpha}{2}}(\overline{Q})}
+\sum_{k=1}^2 |\hat\rho y_{k, d}|_{L^2(Q)}\leq \delta^0,
$$ 
one  can always  find  a leader  control
 $\bar u\in C^{1+\alpha, \frac{1+\alpha}{2}}(\overline{Q})$ so that    $J_1(\cdot, \cdot; \bar u)$  
and  $J_2(\cdot, \cdot; \bar u)$ have a unique  pair of 
 Nash equilibrium $(\bar{v}_1,\bar{v}_2)\in (B_{\rho_1})^2$.  Moreover,   the corresponding  solution 
 $y(\cdot, \cdot; \bar u, \bar{v}_1,\bar{v}_2)$ to $(\ref{a})$
satisfies that $$y(\cdot,T; \bar u, \bar{v}_1,\bar{v}_2)=0\quad \mbox{   in } \Omega.$$
\end{theorem}

In Theorems \ref{th1} and \ref{+6}, 
the locally null controllability problem for the quasi-linear parabolic equation 
(\ref{a}) by hierarchic controls 
is studied. Indeed, 
 it can be reduced to a partial locally  null controllability
for a strongly coupled quasi-linear parabolic system. 
First, for any leader control,  a pair of Nash quasi-equilibrium may be represented as 
two components of solutions  to a coupled  system, in which,  coupling appears in 
principal parts of the equations.  
Then the associated controllability problem is to find a leader control so 
that one component of solutions
to this strongly coupled system reaches the desired target zero for sufficiently small initial value 
$y_0$
and $y_{i, d}$ in the cost functionals.
As for the local  controllability of single quasi-linear parabolic equations,  
the known results 
 usually require the coefficients in principal parts
  to be independent of gradient of solutions, or spacial dimension to be limited (see 
  \cite{MB, f, Liux,  NN, MRN} and the references therein).
In this paper,  we will try to improve  regularity of controls  in  suitable 
H\"older  spaces and 
establish the local controllability results for the strongly coupled quasi-linear parabolic system,
 in which,  
the coefficients in principal parts of  the equations
 contain not only the state itself but also its gradient with general spacial dimension.

The rest of this paper is organized as follows. 
Section 2 is devoted to characterizing  a pair of 
Nash quasi-equilibrium for the functionals $J_1$ and $J_2$ for any leader control.
 In Section 3, we prove an observability inequality for a coupled linear parabolic system.
 The key is to establish a suitable Carleman estimate for the system and  to 
 present  dependence of the estimate on  coefficients of it.  
 In Section 4,  we  improve  regularity of controls and 
 prove a controllability result for a coupled linear parabolic system.
  Sections 5 is devoted to proving the main  results  of this paper, 
 i.e., Theorems \ref{th1} and \ref{+6}.

\section{Characterization of Nash quasi-equilibrium} 

In this section,  as the first step of solving the multi-objective problem, 
 for any given leader control $u\in B_{\rho_1}$,  we prove the  existence of  a pair of 
 Nash quasi-equilibrium 
$(\bar{v}_1,\bar{v}_2)\in (B_{\rho_1})^2$ for
the functionals $J_1$ and $J_2$ 
 in (\ref{d}).    It will be  represented by    a pair of   components of solutions  to 
 a strongly coupled quasi-linear
parabolic system.  

For convenience,   for a differentiable   function $F: \mathbb R\times \mathbb R^n\rightarrow  \mathbb R$,  
we denote by $F_y$ and $\nabla_\zeta F$ the partial derivative and  gradient  
of $F=F(y, \zeta)$   with 
respect to the first and the rest variables,  respectively. 
Further,  we use $\nabla\cdot v$ and  $F_{\zeta_k}$  to represent
the divergence of a vector $v$ and the partial derivative of $F$ with respect   
to the $k+1$-th independent variable for $k=1, \cdots, n$, respectively.

\medskip

The main  result of  this section is given  as follows.
\begin{proposition}\label{DX!}
For any $y_0\in \tilde B_{\rho_1}$ and  leader control $u\in B_{\rho_1}$, 
there exists a $\mu^*>0$, such that for any $\mu_k\geq \mu^*$ $(k=1, 2)$, 
$(\bar v_1, \bar v_2)\in \big(B_{\rho_1}\big)^2$ is a Nash quasi-equilibrium of the 
functionals $J_1$ and $J_2$ in $(\ref{d})$, if and only if 
$$
\bar v_1=\frac{1}{\mu_1} \xi_{1}p_1\quad \mbox{ and } 
\quad
\bar v_2=\frac{1}{\mu_2} \xi_{2}p_2,$$
where $(\bar y, p_1, p_2)\in \big(C^{2+\alpha, 1+\frac{\alpha}{2}}(\overline{Q})\big)^3$ satisfy the following coupled quasi-linear parabolic  system:

\begin{eqnarray}\label{***}
\left\{
\begin{array}{ll}
\bar y_{t}-\sum\limits^n_{i, j=1}\Big(a^{ij}(\bar y,\nabla \bar y)\bar y_{ x_i}\Big)_{x_j}
+f(\bar y, \nabla \bar y)=\xi_{0}u+
\displaystyle\frac{1}{\mu_1}
\xi_{1}^2 p_1
+\frac{1}{\mu_2}
\xi_{2}^2 p_2     &\mbox{ in }   Q,\\
p_{1, t}+\sum\limits^n_{i, j=1} \Big(a^{ij}(\bar y, \nabla \bar y)p_{1, x_i}\Big)_{x_j}
+\sum\limits_{i, j=1}^n \div
\Big[\bar{y}_{x_i} p_{1, x_j}\nabla_{\zeta} a^{i j}(\bar y, \nabla \bar y)\Big]
&\\
\quad-\sum\limits_{i, j=1}^n 
a^{i j}_y(\bar y, \nabla\bar y)\bar{y}_{x_i}p_{1, x_j}
-f_y(\bar y, \nabla\bar y)p_1+\div \Big[p_1\nabla_{\zeta}f(\bar y, \nabla\bar y)\Big]=
\nu_1 \xi_{*}(\bar y-y_{1, d}) &\text{ in }Q,\\
p_{2, t}+\sum\limits^n_{i, j=1} \Big(a^{ij}(\bar y, \nabla \bar y)p_{2, x_i}\Big)_{x_j}
+\sum\limits_{i, j=1}^n \div
\Big[\bar{y}_{x_i} p_{2, x_j}\nabla_{\zeta} a^{i j}(\bar y, \nabla \bar y)\Big]
&\\
\quad-\sum\limits_{i, j=1}^n 
a^{i j}_y(\bar y, \nabla\bar y)\bar{y}_{x_i}p_{2, x_j}
-f_y(\bar y, \nabla\bar y)p_2+\div \Big[p_2\nabla_{\zeta}f(\bar y, \nabla\bar y)\Big]=
\nu_2 \xi_{*}(\bar y-y_{2, d}) &\text{ in }Q,\\

\bar y=p_1=p_2=0  &\text{ on }  \Sigma,\\[2mm]
\bar y(x, 0)=y_0(x),\    p_1(x,  T)=p_2(x, T)=0 &\text{ in }  \Omega.
\end{array}
\right.
\end{eqnarray}
\end{proposition}
\noindent {\bf Proof. }  The whole proof is divided into three parts.

\medskip

{\bf Step 1. }
First, for any leader control $u\in B_{\rho_1}$, 
$v_1,  v_2,  \bar v_1,  \bar v_2\in B_{\rho_1}$ and sufficiently small $r\in (0, 1)$ with 
$\bar v_k+r v_k\in B_{\rho_1}$ for $k=1, 2$, 
let $y_1^r$ and $y_2^r$ be
the solutions to (\ref{a}) associated to  control triples 
$(u, \bar{v}_1+rv_1, \bar{v}_2)$ and $(u, \bar{v}_1, \bar{v}_2+rv_2)$,
respectively.  Then, $y_1^r$ and $y_2^r$ satisfy the following 
quasi-linear parabolic equations: 
\begin{eqnarray*}\label{21}
\left\{
\begin{array}{ll}
y^r_{1, t}-\sum\limits^n_{i,j=1}\Big(a^{ij}(y_1^r,\nabla y_1^r)y^r_{1, x_i}\Big)_{x_j}
+f(y_1^r, \nabla  y_1^r)=\xi_{0}u+\xi_{1}(\bar{v}_1+rv_1)+\xi_{2}\bar{v}_2      &\mbox{ in }   Q,\\
y^r_1=0  &\mbox{ on }  \Sigma,\\[2mm]
y^r_1(x,  0)=y_0(x)  &\mbox{ in }  \Omega,
\end{array}
\right.
\end{eqnarray*}
and
\begin{eqnarray*}\label{22}
\left\{
\begin{array}{ll}
y^r_{2, t}-\sum\limits^n_{i,j=1}\Big(a^{ij}(y_2^r,\nabla y_2^r)y^r_{2, x_i}\Big)_{x_j}
+f(y_2^r, \nabla  y_2^r)=\xi_{0}u+\xi_{1}\bar{v}_1
+\xi_{2}(\bar{v}_2+rv_2)      &\mbox{ in }   Q,\\
y^r_2=0  &\mbox{ on }  \Sigma,\\[2mm]
y^r_2(x,  0)=y_0(x)  &\mbox{ in }  \Omega.
\end{array}
\right.
\end{eqnarray*}
By  well-posedness results of quasi-linear parabolic equations and (\ref{XX1}), 
both $\{y_1^r\}$ and 
$\{y_2^r\}$ are bounded in 
$C^{2+\alpha, 1+\frac{\alpha}{2}}(\overline{Q})$ with respect to  $r$.
Hence,  there exist subsequences of 
$\{y_1^r\}$ and 
$\{y_2^r\}$, which are still denoted by themselves, such that
$$
\lim\limits_{r\rightarrow 0} y^r_1=\lim\limits_{r\rightarrow 0} y^r_2=\bar y
\quad\mbox{ in }C^{2, 1}(\overline{Q}),
$$
where  $\bar y\in C^{2+\alpha, 1+\frac{\alpha}{2}}(\overline{Q})$ is the solution to (\ref{a}) corresponding to  control 
triple $(u, \bar  v_1,   \bar  v_2)$.

 Next, set 
$$
y_1(\cdot, \cdot)=\lim\limits_{r\rightarrow 0}\frac{y_1^r(\cdot, \cdot)-\bar y(\cdot, \cdot)}{r}
\  \mbox{ and }\  
y_2(\cdot, \cdot)=\lim\limits_{r\rightarrow 0}\frac{y_2^r(\cdot, \cdot)-\bar y(\cdot, \cdot)}{r}.$$
Then,  $y_1$ and      $y_2$ satisfy   the following 
partial differential  equations,  respectively:
\begin{eqnarray}\label{23}
\left\{
\begin{array}{ll}
y_{1, t}-\sum\limits^n_{i, j=1} \Big(a^{ij}(\bar y, \nabla \bar y)y_{1, x_i}\Big)_{x_j}
-\sum\limits_{i, j=1}^n \Big(\bar{y}_{x_i} \nabla_\zeta a^{i j}(\bar{y}, \nabla\bar y)\cdot 
\nabla y_1\Big)_{x_j}&\\
\quad-\sum\limits_{i, j=1}^n \Big(a^{i j}_y(\bar y, \nabla\bar y)\bar{y}_{x_i}y_1\Big)_{x_j}
+f_y(\bar y, \nabla\bar y)y_1+\nabla_{\zeta}f(\bar y, \nabla\bar y)\cdot\nabla y_1=
\xi_{1}v_1  &\mbox{ in }Q,\\
y_1=0  &\mbox{ on }  \Sigma,\\[2mm]
y_1(x,  0)=0 &\mbox{ in }  \Omega,
\end{array}
\right.
\end{eqnarray}
and
\begin{eqnarray}\label{24}
\left\{
\begin{array}{ll}
y_{2, t}-\sum\limits^n_{i, j=1} \Big(a^{ij}(\bar y, \nabla \bar y)y_{2, x_i}\Big)_{x_j}
-\sum\limits_{i, j=1}^n \Big(\bar{y}_{x_i} \nabla_\zeta a^{i j}(\bar{y}, \nabla\bar y)\cdot 
\nabla y_2\Big)_{x_j}&\\
\quad-\sum\limits_{i, j=1}^n \Big(a^{i j}_y(\bar y, \nabla\bar y)\bar{y}_{x_i}y_2\Big)_{x_j}
+f_y(\bar y, \nabla\bar y)y_2+\nabla_{\zeta}f(\bar y, \nabla\bar y)\cdot\nabla y_2=
\xi_{2}v_2  &\mbox{ in }Q,\\
y_2=0  &\mbox{ on }  \Sigma,\\[2mm]
y_2(x,  0)=0 &\mbox{ in }  \Omega.
\end{array}
\right.
\end{eqnarray}

Note that  in the first equation of (\ref{23}),
\begin{eqnarray*}
&&\sum\limits_{i, j=1}^n \Big(\bar{y}_{x_i} \nabla_\zeta a^{i j}(\bar{y}, 
\nabla\bar y)\cdot 
\nabla y_1\Big)_{x_j}
=\sum\limits_{i, j, \ell=1}^n \Big(\bar{y}_{x_i} a^{i j}_{\zeta_{\ell}}(\bar{y}, \nabla\bar y)
y_{1, x_{\ell}}\Big)_{x_j}\\
&&=\sum\limits_{i, j=1}^n \Big(\sum\limits_{\ell=1}^n
\bar{y}_{x_{\ell}} a^{\ell j}_{\zeta_i}(\bar{y}, \nabla\bar y)
y_{1, x_i}\Big)_{x_j}.
\end{eqnarray*}
By the previous notations,  $a^{i j}_{\zeta_{\ell}}$ denotes the partial derivative of $a^{i j}=a^{i j}(y, 
\zeta_1, \cdots, \zeta_n)$ with respect to the   variable $\zeta_{\ell}$. Set 
\begin{equation}\label{DX5}
c^{i j}=\displaystyle\frac{1}{2}\sum\limits_{\ell=1}^n
\bar{y}_{x_{\ell}}[a^{\ell j}_{\zeta_i}(\bar{y}, \nabla\bar y)
+a^{\ell i}_{\zeta_j}(\bar{y}, \nabla\bar y)].
\end{equation}
Then $c^{i j}=c^{j i}$ in $\overline{Q}$ for $i, j=1, 2, \cdots, n$. 
Since $a^{i j}\in C^4(\mathbb R\times\mathbb R^n)$ and 
$\bar y\in C^{2+\alpha, 1+\frac{\alpha}{2}}(\overline{Q})$, 
we have that
$c^{i j}\in \widetilde C_1^{\alpha, 
\frac{\alpha}{2}}(\overline{Q})$
and $$
\sum\limits_{i, j=1}^n \Big(\bar{y}_{x_i} \nabla_\zeta a^{i j}(\bar{y}, 
\nabla\bar y)\cdot 
\nabla y_1\Big)_{x_j}=\sum\limits_{i, j=1}^{n}(c^{i j} y_{1, x_i})_{x_j}.$$
By (\ref{XX1}), 
$
|\bar y|_{C^{2+\alpha, 1+\frac{\alpha}{2}}(\overline{Q})}
\leq C(n, T, \Omega, \alpha, a^{i j}, f)\rho_1.$
 Hence, 
  when $\rho_1$ is small enough, 
$|\bar y|_{C^{2+\alpha, 1+\frac{\alpha}{2}}(\overline{Q})}$ is sufficiently small.
Without loss of generality,     assume that
\begin{equation}\label{DX2}
\Big|\sum\limits_{i, j=1}^{n} c^{i j}(x, t)\zeta_i \zeta_j\Big|\leq 
\frac{\rho_0}{2}|\zeta|^2, \quad 
\forall\ (x, t, \zeta)=(x, t, \zeta_1, \zeta_2, \cdots, \zeta_n)\in \overline{Q}\times\mathbb R^n.
\end{equation}
This implies that (\ref{23})
is a quasi-linear parabolic equation and  its solution 
$y_1\in C^{2+\alpha, 1+\frac{\alpha}{2}}(\overline{Q})$. Similarly, 
(\ref{24})
is also a quasi-linear parabolic equation and  its solution
$y_2\in C^{2+\alpha, 1+\frac{\alpha}{2}}(\overline{Q})$.

\medskip

{\bf Step 2. } By Definition \ref{+1},  
$(\bar{v}_1,\bar{v}_2)\in  (B_{\rho_1})^2$  is the 
Nash quasi-equilibrium of $J_1$  and  $J_2$, 
 if and only if  
\begin{eqnarray}\label{X1}
\begin{array}{ll}
&\displaystyle J_{1, v_1}(\bar{v}_1,\bar{v}_2; u)v_1
=\lim\limits_{r\rightarrow 0} \frac{J_1(\bar{v}_1+rv_1, \bar{v}_2; u)-
J_1(\bar v_1, \bar v_2; u)}{r}\\[3mm]
&\displaystyle=
\lim\limits_{r\rightarrow0}\Big[
\frac{\mu_1}{2}
\int^T_0\int_{\omega_1}
\frac{(\bar{v}_1+rv_1)^2-\bar{v}_1^2}{r}dxdt
+
\frac{\nu_1}{2}
\int_Q \xi_{*}
\frac{(y^r_1-y_{1,d})^2-(\bar{y}-y_{1,d})^2}{r}dxdt\Big]\\[3mm]
&\displaystyle=\mu_1 \int^T_0\int_{\omega_1}
\bar{v}_1 v_1dxdt
+\nu_1
\int_Q \xi_{*} (\bar y-y_{1,d})y_1dxdt=0,
\end{array}
\end{eqnarray}
where $y_1$ is the solution to (\ref{23}). 
Similarly, 
\begin{equation}\label{X2}
J_{2, v_2}(\bar{v}_1,\bar{v}_2; u)v_2=
\mu_2 \int^T_0\int_{\omega_2}
\bar{v}_2 v_2dxdt
+\nu_2
\int_Q \xi_{*} (\bar y-y_{2,d})y_2dxdt=0, 
\end{equation}
where $y_2$ is the solution to (\ref{24}).

In order to express the Nash quasi-equilibrium $(\bar v_1, \bar v_2)$, for  
 the solution $\bar y$  to (\ref{a}) corresponding to  control 
triple $(u, \bar  v_1,   \bar  v_2)$,  introduce
 the following linear partial differential equations:
\begin{eqnarray}\label{25}
\left\{
\begin{array}{ll}
p_{1, t}+\sum\limits^n_{i, j=1} \Big(a^{ij}(\bar y, \nabla \bar y)p_{1, x_i}\Big)_{x_j}
+\sum\limits_{i, j=1}^n \div
\Big[\bar{y}_{x_i} p_{1, x_j}\nabla_{\zeta} a^{i j}(\bar y, \nabla \bar y)\Big]
&\\
\quad-\sum\limits_{i, j=1}^n 
a^{i j}_y(\bar y, \nabla\bar y)\bar{y}_{x_i}p_{1, x_j}
-f_y(\bar y, \nabla\bar y)p_1+\div \Big[p_1\nabla_{\zeta}f(\bar y, \nabla\bar y)\Big]=
\nu_1 \xi_{*}(\bar y-y_{1, d}) &\mbox{ in }Q,\\
p_1=0  &\mbox{ on }  \Sigma,\\[2mm]
p_1(x,  T)=0 &\mbox{ in }  \Omega,
\end{array}
\right.
\end{eqnarray}
and
\begin{eqnarray}\label{26}
\left\{
\begin{array}{ll}
p_{2, t}+\sum\limits^n_{i, j=1} \Big(a^{ij}(\bar y, \nabla \bar y)p_{2, x_i}\Big)_{x_j}
+\sum\limits_{i, j=1}^n \div
\Big[\bar{y}_{x_i} p_{2, x_j}\nabla_{\zeta} a^{i j}(\bar y, \nabla \bar y)\Big]
&\\
\quad-\sum\limits_{i, j=1}^n 
a^{i j}_y(\bar y, \nabla\bar y)\bar{y}_{x_i}p_{2, x_j}
-f_y(\bar y, \nabla\bar y)p_2+\div \Big[p_2\nabla_{\zeta}f(\bar y, \nabla\bar y)\Big]=
\nu_2 \xi_{*}(\bar y-y_{2, d}) &\mbox{ in }Q,\\
p_2=0  &\mbox{ on }  \Sigma,\\[2mm]
p_2(x,  T)=0 &\mbox{ in }  \Omega.
\end{array}
\right.
\end{eqnarray}

Note that  for $k=1, 2$, by (\ref{DX5}),
\begin{eqnarray*}
&&\sum\limits_{i, j=1}^n \div
\Big[\bar{y}_{x_i} p_{k, x_j}\nabla_{\zeta} a^{i j}(\bar y, \nabla \bar y)\Big]
=\sum\limits_{i, j, \ell=1}^{n} \big(\bar y_{x_\ell} p_{k, x_j} a^{\ell j}_{\zeta_i}(\bar y, \nabla \bar y
)\big)_{x_i}=\sum\limits_{i, j=1}^{n}(c^{i j}p_{k, x_j})_{x_i}.
\end{eqnarray*}
By (\ref{DX2}),  both (\ref{25}) and (\ref{26}) are  linear 
parabolic equations of second order with respect to $p_1$ and $p_2$, respectively. Moreover, 
by the conditions for $\xi_*$  in (\ref{d}), 
 compatibility conditions hold.
Hence, (\ref{25}) and (\ref{26}) have a unique solution $p_1$ and $p_2$ in 
$C^{2+\alpha, 1+\frac{\alpha}{2}}(\overline{Q})$, respectively.
By the Schauder estimate for  linear 
parabolic equations of second order, 
\begin{eqnarray}\label{DX3}
\begin{array}{rl}
&\displaystyle |p_1|_{C^{2+\alpha, 1+\frac{\alpha}{2}}(\overline{Q})}
+\displaystyle |p_2|_{C^{2+\alpha, 1+\frac{\alpha}{2}}(\overline{Q})}\\[2mm]
&
\leq C(T, \Omega, n, \rho_1, \nu_1, \nu_2, \alpha, a^{i j}, f)\big(|\bar y|_{C^{\alpha, 
\frac{\alpha}{2}}(\overline{Q})}+
|y_{1, d}|_{C^{\alpha, 
\frac{\alpha}{2}}(\overline{Q})}
+
|y_{2, d}|_{C^{\alpha, 
\frac{\alpha}{2}}(\overline{Q})}
\big)\\[2mm]
&\displaystyle \leq C(T, \Omega, n, \rho_1, \nu_1, \nu_2, \alpha, 
a^{i j}, f,  y_{1, d}, y_{2, d})\deq C_1.
\end{array}
\end{eqnarray}

By  the duality between (\ref{23}) and  (\ref{25}),  it holds that
\begin{equation}\label{YYY}
\displaystyle\int_Q
\xi_{1} p_1 v_1 dxdt
=-\int_Q \nu_1\xi_{*}(\bar y-y_{1, d})y_1dxdt.
\end{equation}
Similarly,  by  the duality between   (\ref{24})  and (\ref{26}),  it follows that
\begin{equation}\label{YYY*}
\displaystyle\int_Q
\xi_{2} p_2 v_2 dxdt=-\int_Q \nu_2\xi_{*}(\bar y-y_{2, d})y_2dxdt.
\end{equation}

Combining the above two equalities with (\ref{X1})  and (\ref{X2}),  we obtain that  
$$
\displaystyle\int_Q v_k(\xi_{k}p_k-\mu_k \bar v_k \mathbb I_{\omega_k})dxdt=0, 
\qquad\forall\ v_k\in B_{\rho_1} \mbox{ and } k=1, 2,
$$
where $\mathbb I_{\omega_k}$ denotes the characteristic functions on the set 
$\omega_k$.  Hence,
$
\bar v_k=\frac{1}{\mu_k} \xi_{k}p_k,$ $\forall\ k=1, 2,
$
where $p_1$ and $p_2$ are the solutions to (\ref{25}) and (\ref{26}). 
For the constant $C_1$ defined in (\ref{DX3}),
set 
\begin{equation}\label{XX2}
\mu^*=2C_1/\rho_1,
\end{equation} for the constant $C_1$ defined in (\ref{DX3}). Then, 
for any $\mu_1, \mu_2\geq\mu^*$,  we have $\bar v_1, \bar v_2\in B_{\rho_1}$.
Further,   
$(\bar y,   p_1,  p_2)\in (C^{2+\alpha, 1+\frac{\alpha}{2}}(\overline{Q}))^3$ are solutions to the coupled system (\ref{***}).

\medskip

{\bf Step 3. } On the other hand,  for any $\mu_1, \mu_2\geq\mu^*=2C_1/\rho_1$, set 
\begin{equation}\label{DX}
\bar v_1=\frac{1}{\mu_1} \xi_{1}p_1\quad \mbox{ and } 
\quad
\bar v_2=\frac{1}{\mu_2} \xi_{2}p_2,
\end{equation} where 
$p_1$ and $p_2$ are  components  of solutions  to (\ref{***}).

 By the duality between (\ref{23})  and  (\ref{25}), and between  (\ref{24})  and (\ref{26}),  
we have that  both (\ref{YYY}) and  (\ref{YYY*}) hold. 
Furthermore,  by (\ref{DX}), 
they imply that both (\ref{X1})
and (\ref{X2}) are true,  that is,  $(\bar  v_1, \bar v_2)\in (B_{\rho_1})^2$ defined in (\ref{DX}) is  
 the Nash quasi-equilibrium
of $J_1$ and $J_2$. This completes the proof of Proposition \ref{DX!}.
\endpf

\medskip

By Proposition \ref{DX!},  for any given leader control $u\in B_{\rho_1}$, 
there exists a unique Nash quasi-equilibrium $(\bar v_1, \bar v_2)\in (B_{\rho_1})^2$ (as 
a pair of follower controls) and 
it can be represented as the solution components for a strongly coupled
 quasi-linear parabolic system (\ref{***}) with sufficiently large $\mu_1$ and $\mu_2$.

 For a more concise expression,  we introduce the following notations.
   For  $i, j=1, 2,  \cdots, n$,   set 
\begin{eqnarray}\label{dl1}
\begin{array}{ll}
c^{i j}(\bar y, \nabla\bar  y)=\displaystyle\frac{1}{2}\sum\limits_{\ell=1}^n 
\bar{y}_{x_\ell}  [a^{\ell i}_{\zeta_j}(\bar y, \nabla \bar y)+
a^{\ell j}_{\zeta_i}(\bar y, \nabla \bar y)] \  (\mbox{as in } (\ref{DX5})),&\\[3mm]
e^j(\bar y, \nabla\bar  y)=\displaystyle-\sum\limits_{i=1}^n 
a^{i j}_{y}(\bar y, \nabla \bar y)\bar y_{x_i}+f_{\zeta_j}(\bar y,  \nabla\bar y),&\\[4mm]  
d_0(\bar y, \nabla\bar  y)=-f_y(\bar y, \nabla\bar y)+
\sum\limits_{\ell=1}^n\big[f_{\zeta_{\ell}}(\bar y, \nabla\bar  y)\big]_{x_\ell},&\\[4mm]
\mbox{and } A^{i j}(\bar y, \nabla\bar  y)=a^{ij}(\bar y, \nabla \bar y)+
c^{i j}(\bar y, \nabla \bar y).
\end{array}
\end{eqnarray}
Then, the coupled system (\ref{***}) may be rewritten as follows:
\begin{eqnarray}\label{dl2}
\left\{
\begin{array}{ll}
\bar y_{t}-\sum\limits^n_{i, j=1}\big(a^{ij}(\bar y,\nabla \bar y)\bar y_{ x_i}\big)_{x_j}
+f(\bar y, \nabla \bar y)=\xi_{0}u+
\displaystyle\frac{1}{\mu_1}
\xi_{1}^2 p_1
+\frac{1}{\mu_2}
\xi_{2}^2 p_2     &\mbox{ in }   Q,\\
p_{1, t}+\sum\limits^n_{i, j=1} \big(A^{ij}(\bar y, \nabla \bar y)
p_{1, x_i}\big)_{x_j} &\\
\quad\quad\quad\quad\quad+\sum\limits_{j=1}^n 
e^j(\bar y, \nabla \bar y)  p_{1, x_j}+d_0(\bar y, \nabla \bar y)
p_1=
\nu_1 \xi_{*}(\bar y-y_{1, d}) &\mbox{ in }Q,\\
p_{2, t}+\sum\limits^n_{i, j=1} \Big(A^{ij}(\bar y, \nabla \bar y)
p_{2, x_i}\Big)_{x_j}&\\
\quad\quad\quad\quad\quad+\sum\limits_{j=1}^n 
e^j(\bar y, \nabla \bar y)  p_{2, x_j}+d_0(\bar y, \nabla \bar y)
p_2=
\nu_2 \xi_{*}(\bar y-y_{2, d}) &\mbox{ in }Q,\\
\bar y=p_1=p_2=0  &\mbox{ on }  \Sigma,\\[2mm]
\bar y(x, 0)=y_0(x), \ p_1(x,  T)=p_2(x, T)=0 &\mbox{ in }  \Omega.
\end{array}
\right.
\end{eqnarray}

In (\ref{dl2}), the coefficients in principal parts of the second and third equations 
contain the solution component $\bar y$. Hence, it is a strongly coupled 
quasi-linear parabolic system.
The aim of  this paper is to find a leader  control $u\in B_{\rho_1}$, 
such that the solution triple $(\bar  y, p_1, p_2)$ to the coupled system
 (\ref{dl2}) satisfies that 
$
\bar y(x, T)=0 \mbox{ in }\Omega.$
To solve this   controllability problem for the coupled  quasi-linear parabolic system, 
we shall  utilize the fixed point technique and first study  the associated  controllability problem
 for a coupled linear parabolic system
  in the frame of classical solutions. 
As preliminaries,  the next section is devoted to 
establishing an observability problem for a coupled linear parabolic system.

\section{Observability inequality  for coupled linear parabolic systems}

In order to establish the desired controllability result for the linearized system of (\ref{dl2}), 
in this section,  we first derive a global Carleman estimate for the following coupled  linear parabolic system:
\begin{equation}\label{s}
  \begin{cases}
  \begin{split}
&\varphi_t+\sum\limits^n_{i,j=1}(b^{ij}\varphi_{x_i})_{x_j}+
\sum\limits_{j=1}^{n} (f_j\varphi)_{x_j}-f_0\varphi+
 \xi_{*}(\nu_1\theta_1+\nu_2\theta_2)=0  &\mbox{ in }Q,&\\
&\theta_{1, t}-\sum\limits^n_{i,j=1}(B^{ij}\theta_{1, x_i})_{x_j}+
\sum\limits_{j=1}^{n} (g_j\theta_1)_{x_j}-g_0\theta_1+
 \frac{1}{\mu_1}\xi_{1}^2\varphi=0   &\mbox{ in }  Q,&\\
&\theta_{2, t}-\sum\limits^n_{i,j=1}(B^{ij}\theta_{2, x_i})_{x_j}+
\sum\limits_{j=1}^{n} (g_j\theta_2)_{x_j}-g_0\theta_2+
 \frac{1}{\mu_2}\xi_{2}^2\varphi=0    &\mbox{ in }  Q,&\\
&\varphi=\theta_1= \theta_2=0   &\mbox{ on } \Sigma,&\\[2mm]
&\varphi(x, T)=\varphi_{T}(x), \  \theta_1(x, 0)=\theta_2(x, 0)=0 &\mbox{ in } \Omega,&
  \end{split}
  \end{cases}
\end{equation}
where   
$b^{i j}, B^{i j}\in C^{1, 1}(\overline{Q})\cap \widetilde 
C_1^{\alpha, \frac{\alpha}{2}}(\overline{Q})$  with 
$b^{i j}(x, t)=b^{j i}(x, t)$  and $B^{i j}(x, t)=B^{j i}(x, t)$ in $\overline{Q}$, 
$f_j, g_j\in \widetilde C_1^{\alpha, \frac{\alpha}{2}}(\overline{Q})$ for $i, j=1, \cdots, n$, 
 $f_0, g_0\in C^{\alpha, \frac{\alpha}{2}}(\overline{Q})$, 
 and $\varphi_T\in C^{2+\alpha}_0(\Omega)$.   Moreover,   
for the positive  constant $\rho_0$ and any $(x, t,  \zeta)=(x, t,  \zeta_1, \cdots, \zeta_n)\in
 \overline{Q}\times\mathbb R^n$,  
$$
\sum\limits_{i, j=1}^n  b^{i j}(x, t)\zeta_i\zeta_j\geq \rho_0  |\zeta|^2 \quad\mbox{and} \quad
\sum\limits_{i, j=1}^n  B^{i j}(x, t)\zeta_i\zeta_j\geq \rho_0  |\zeta|^2.
$$

	Set  $\psi=\nu_1 \theta_1+\nu_2\theta_2$.  Then (\ref{s}) can be  simplified  as follows:
\begin{equation}\label{t}
  \begin{cases}
  \begin{split}
&\varphi_t+\sum\limits^n_{i,j=1}(b^{ij}\varphi_{x_i})_{x_j}+
\sum\limits_{j=1}^{n} (f_j\varphi)_{x_j}-f_0\varphi+
\xi_{*}\psi=0  &\mbox{ in }Q,&\\
&\psi_{t}-\sum\limits^n_{i,j=1}(B^{ij}\psi_{x_i})_{x_j}+
\sum\limits_{j=1}^{n} (g_j\psi)_{x_j}-g_0\psi+\big(
\frac{\nu_1}{\mu_1}\xi_{1}^2+\frac{\nu_2}{\mu_2}\xi_{2}^2\big)\varphi=0   &\mbox{ in }  Q,&\\
&\varphi=\psi=0 \quad  &\mbox{ on } \Sigma,&\\
& \varphi(x, T)=\varphi_{T}(x),\   \psi(x, 0)=0 \quad    &\mbox{ in } \Omega.&
  \end{split}
  \end{cases}
\end{equation}

First, we recall a global  Carleman estimate for the following linear parabolic equation:
\begin{equation}\label{u}
  \begin{cases}
  \begin{split}
  & v_t+\sum\limits_{i,j=1}^n (b^{ij}v_{x_i} )_{x_j}+\sum\limits_{j=1}^n (f_j v)_{x_j}+f_0 v=h_0 \qquad &\mbox{ in }  Q,& \\
  & v=0 \qquad &\mbox{ on } \Sigma,&\\
  & v(x, T)=v_T(x) \qquad  &\mbox{ in } \Omega,& \\
  \end{split}
  \end{cases}
\end{equation}
where $h_0 \in L^2(Q)$ and $v_T \in L^2(\Omega)$.
Let $\mathcal{O}$  be any given nonempty open  subset of $\Omega$
and $\eta\in C^4(\overline{\Omega})$ satisfy that 
$$
\eta>0 \mbox{ in } \Omega,\quad \eta=0 \mbox{ on }\partial\Omega\ \   \mbox{ and }\  \   
|\nabla\eta|>0 \mbox{ in } \overline{\Omega\setminus\mathcal{O}}.$$
For any  parameter  $\mu\geq 1$,   introduce the following functions, which will be used in Carleman 
estimates:
\begin{eqnarray*}
	&&\beta(x,t)=\frac{e^{\mu\eta(x)}}{t(T-t)},\quad \beta_0(t)=\frac{1}{t(T-t)},\\
&&\nu(x,t)=\big(e^{\mu\eta(x)}-
e^{2\mu|\eta|_{C(\overline{\Omega})}}\big)\beta_0(t)
\mbox{ and }
 \nu_0(t)=\big(1-e^{2\mu|\eta|_{C(\overline{\Omega})}}\big) \beta_0(t),
\end{eqnarray*}
and write
$$B_0=\sum\limits^n_{i,j=1}(1+|b^{ij}|^2_{C^{1,1}(\overline{Q})})\  \mbox{ and }\  
D_0=\sum\limits^n_{j=1}|f_j|^2_{C^{1,0}(\overline{Q})}+|f_0|^2_{C(\overline Q)}.$$
Then,  a known Carleman estimate for  (\ref{u})  is given as follows (see \cite[Proposition 3.2]{Liu}).
\begin{lemma}\label{prop1}
There  exists a positive constant $C$,  such that for any  $\mu\geq C B_0$ and 
$\lambda \geq  C(D_0+e^{2\mu|\eta|_{C(\overline{\Omega})}})$,
 any solution to the equation $(\ref{u})$ satisfies
\begin{equation}\label{v}
  \begin{split}
&\int_Q e^{2\lambda\nu}\Big(\lambda\mu^2\beta|\nabla v|^2+\lambda^3\mu^4\beta^3v^2\Big)
dxdt
\\
&\leq  C\Big(1+\sum\limits_{i,j=1}^n|b^{ij}|^2_{C(\overline{Q})}+\sum\limits_{j=1}^n|f_j|^2_{C(\overline{Q})}+|f_0|_{C(\overline Q)}\Big)
\int_0^T\int_{\mathcal{O}} e^{2\lambda\nu}\lambda^3\mu^4\beta^3v^2dxdt\\
& \q +C\int_Q e^{2\lambda\nu}h_0^2dxdt.
\end{split}
\end{equation}
\end{lemma}

Based on Lemma \ref{prop1}, using similar arguments to \cite{fur}, one can immediately get the following Carleman estimate. 
\begin{lemma}\label{lem1}
There  exists a positive constant $C$,  such that for any  $\mu\geq C B_0$ and $\lambda 
\geq  C(D_0+e^{2\mu|\eta|_{C(\overline{\Omega})}})$,
any solution to the equation $(\ref{u})$ satisfies
\begin{equation}\label{DX8}
  \begin{split}
&\int_Q e^{2\lambda\nu}(\lambda \mu^2\beta|\nabla v|^2+\lambda^3\mu^4\beta^3v^2)dxdt
+\int_Q e^{2\lambda\nu} \frac{1}{\lambda\beta} \Big[v_t^2  +\sum\limits^n_{i,j=1}|v_{x_ix_j}|^2\Big]dxdt
\\
&\leq C\Big(1+\sum\limits_{i,j=1}^n |b^{ij}|^8_{C^{0, 1}(\overline{Q})}
 +\sum\limits_{j=1}^n |f_j|^8_{C^{1, 0}(\overline{Q})}+|f_0|^6_{C(\overline{Q})}\Big)\\
& \q \q\q\cdot\Big(\int_0^T\int_{\mathcal{O}} e^{2\lambda\nu}\lambda^3\mu^4\beta^3v^2dxdt
+\int_Q e^{2\lambda \nu}h_0^2 dxdt\Big).
\end{split}
\end{equation}
\end{lemma}
For simplicity of notations,  we assume that
$$\sum\limits_{i,j=1}^n\!\! \big(|b^{ij}|_{C^{1, 1}(\overline{Q})}\!+\!|B^{ij}|_{C^{1, 1}(\overline{Q})}\big)
\!\!+\!\!\sum\limits_{j=1}^n\! (|f_j|_{C^{1, 0}(\overline{Q})}\!+\!|g_j|_{C^{1, 0}(\overline{Q})})\!+\!|f_0|_{C(\overline{Q})}
\!+\!|g_0|_{C(\overline{Q})}\leq \mathcal{L}
$$
and $
B=\sum\limits^n_{i,j=1}(1+|b^{ij}|^2_{C^{1,1}(\overline{Q})}+|B^{ij}|^2_{C^{1,1}(\overline{Q})})$, 
and  use $C(\mathcal{L})$ to denote  positive constants, depending on $n,T,  \Omega, \mathcal{O}, \nu_1, \nu_2$
and $\mathcal{L}$, which may be different in different places.

\smallskip

Next, we apply Lemma \ref{lem1} to derive a Carleman estimate for the coupled system (\ref{t}).

\begin{proposition}\label{prop2}
Assume that  $\tilde\omega_0\cap \omega'\neq\emptyset$.  Then there exists a positive constant $C$,  such that for any $\mu\geq CB$ and $\lambda
\geq C(\mathcal{L})e^{2\mu|\eta|_{C(\overline{\Omega})}}$,
 any solution pair $(\varphi, \psi)$ to  the coupled system $(\ref{t})$ satisfies
\begin{equation}\label{w}
  \begin{split}
  &\int_Q e^{2\lambda\nu} (\lambda \mu^2\beta|\nabla\varphi|^2
  +\lambda^3\mu^4\beta^3\varphi^2) dxdt
+ \int_Q e^{2\lambda\nu} \frac{1}{\lambda\beta} 
\Big[\varphi_t^2  +\sum\limits^n_{i,j=1}|\varphi_{x_ix_j}|^2\Big]dxdt \\
&\quad+\int_Q e^{2\lambda\mu} (\lambda\mu^2\beta|\nabla\psi|^2
+\lambda^3\mu^4\beta^3\psi^2) dxdt
 +\int_Q e^{2\lambda\nu} \frac{1}{\lambda\beta} \Big[\psi_t^2
+\sum\limits^n_{i,j=1}|\psi_{x_ix_j}|^2\Big] dxdt \\
&\leq C(\mathcal{L})\int_0^T \int_{\tilde\omega_0}e^{2\lambda\nu}\lambda^7\mu^8\beta^7\varphi^2   dxdt.
\end{split}
\end{equation}
\end{proposition}

\noindent {\bf Proof.} First, 
we apply the  estimate $(\ref{DX8})$ to the first and second equation of $(\ref{t})$, respectively, for $\mathcal{O}=\omega''$, where $\omega''$ is a nonempty open subset of 
$\tilde\omega_0\cap \omega'$ such that  $\overline{\omega''}\subseteq
\tilde\omega_0\cap \omega'$. It follows that 
\begin{eqnarray*}
&&\int_Q e^{2\lambda\nu}(\lambda\mu^2\beta|\nabla \varphi|^2+\lambda^3\mu^4\beta^3\varphi^2)dxdt
+\int_Q e^{2\lambda\nu}\frac{1}{\lambda\beta}\Big[\varphi_t^2
+\sum\limits^n_{i,j=1}|\varphi_{x_ix_j}|^2\Big]  
 dxdt \\
&&\leq C(\mathcal{L}) \Big(
\int_0^T\int_{\omega''} e^{2\lambda\nu}\lambda^3\mu^4\beta^3\varphi^2 dxdt
+\int_Q e^{2\lambda\nu}|\psi\xi_{*}|^2 dxdt\Big),
\end{eqnarray*}
and
\begin{eqnarray*}
&&\int_Q e^{2\lambda\nu}(\lambda\mu^2\beta|\nabla \psi|^2+\lambda^3\mu^4\beta^3\psi^2)dxdt
 +\int_Q e^{2\lambda\nu}\frac{1}{\lambda\beta}
 \Big[\psi_t^2+\sum\limits^n_{i,j=1}|\psi_{x_ix_j}|^2\Big] dxdt\\
&&\leq C(\mathcal{L})\Big(
\int_0^T\int_{\omega''} e^{2\lambda\nu}\lambda^3\mu^4\beta^3\psi^2 dxdt
+\int_Q e^{2\lambda\nu}\Big|\frac{\nu_1}{\mu_1}\xi^2_{1}
+\frac{\nu_2}{\mu_2}\xi^2_{2}\Big|^2\varphi^2 dxdt\Big).
\end{eqnarray*}
We add the above two estimates together.  When $\lambda\geq C(\mathcal{L})$, 
 the last term in two estimates can be absorbed by 
the left terms and we obtain that
\begin{eqnarray}\label{DX9}
\begin{array}{rl}
&\displaystyle\int_Q e^{2\lambda\nu}(\lambda\mu^2\beta|\nabla \varphi|^2+\lambda^3\mu^4\beta^3\varphi^2)dxdt
+\int_Q e^{2\lambda\nu}\frac{1}{\lambda\beta}\Big[\varphi_t^2
+\sum\limits^n_{i,j=1}|\varphi_{x_ix_j}|^2\Big]  
 dxdt \\
 &\displaystyle+\int_Q e^{2\lambda\nu}(\lambda\mu^2\beta|\nabla \psi|^2+\lambda^3\mu^4\beta^3\psi^2)dxdt
 +\int_Q e^{2\lambda\nu}\frac{1}{\lambda\beta}
 \Big[\psi_t^2+\sum\limits^n_{i,j=1}|\psi_{x_ix_j}|^2\Big] dxdt\\
 &\displaystyle\leq C(\mathcal{L})
\Big(
\int_0^T\int_{\omega''} e^{2\lambda\nu}\lambda^3\mu^4\beta^3\varphi^2 dxdt+
\int_0^T\int_{\omega''} e^{2\lambda\nu}\lambda^3\mu^4\beta^3\psi^2 dxdt\Big).
\end{array}
\end{eqnarray}

Next,  we estimate the last term in (\ref{DX9}). To this aim, 
choose a function $\eta_1 \in C_0^\infty(\tilde\omega_0\cap \omega')$
 satisfying that  $\eta_1=1$ in $\omega''$  and $\eta_{1, x_i}/\eta_1^{\frac{1}{2}}, 
 \eta_{1, x_i x_j}/\eta_1^{\frac{1}{2}} \in L^\infty(\Omega)$ for $i, j=1, \cdots,  n$.
We multiply the first equation of $(\ref{t})$ by $e^{2\lambda\nu}\lambda^3\mu^4\beta^3\eta_1\psi$ and get that
\begin{equation}
  \begin{split}
  \nonumber
&\int_Q e^{2\lambda\nu} \lambda^3\mu^4\beta^3\eta_1\psi^2  dxdt \\
&\displaystyle=
\int_Q  
\Big[-\varphi_t-\sum\limits^n_{i,j=1}(b^{ij}\varphi_{x_i})_{x_j}-\sum\limits^n_{j=1}(f_{j}\varphi)_{x_j}+f_0\varphi 
\Big]e^{2\lambda\nu}\lambda^3\mu^4\beta^3\eta_1\psi dxdt \\
&=\int_Q  e^{2\lambda\nu}\lambda^3\mu^4\beta^3\eta_1f_0\varphi\psi dxdt
+\int_Q e^{2\lambda\nu}\lambda^3\mu^4\beta^3\eta_1\varphi\psi_t dxdt
+\int_Q \varphi\psi\eta_1(e^{2\lambda\nu}\lambda^3\mu^4\beta^3)_t dxdt \\
&\quad+\int_Q \sum\limits^n_{j=1} f_j\psi_{x_j} e^{2\lambda\nu}\lambda^3\mu^4
\beta^3\eta_1\varphi dxdt
+\int_Q\sum\limits^n_{j=1} f_j(e^{2\lambda\nu}\lambda^3\mu^4\beta^3\eta_1)_{x_j} \varphi\psi dxdt\\
& \quad-\sum\limits^n_{i,j=1}\int_Q b^{ij}\varphi\psi(e^{2\lambda\nu}\lambda^3\mu^4\beta^3\eta_1)_{x_ix_j}dxdt
-\sum\limits^n_{i,j=1}\int_Q (b^{ij}\psi)_{x_i}
(e^{2\lambda\nu}\lambda^3\mu^4\beta^3\eta_1)_{x_j}\varphi dxdt \\
& \quad-\sum\limits^n_{i,j=1}\int_Q  (b^{ij} \psi_{x_i})_{x_j}
 \varphi e^{2\lambda\nu}\lambda^3\mu^4\beta^3\eta_1 dxdt
-\sum\limits^n_{i,j=1}\int_Q b^{ij}\psi_{x_j}\varphi(e^{2\lambda\nu}\lambda^3\mu^4\beta^3\eta_1)_{x_i} dxdt \\
&=\sum\limits^9_{j=1}I_j.
\end{split}
\end{equation}

   In the following, we estimate $I_j$ in the above equality, respectively, for $j=1, 2, \cdots, 9$.  
   For any $\varepsilon>0$ and $\lambda\geq C e^{2\mu|\eta|_{C(\overline{\Omega})}}$, 
\begin{eqnarray*}
&&|I_1|=\Big|\int_Q  e^{2\lambda\nu}\lambda^3\mu^4\beta^3\eta_1 f_0\varphi\psi dxdt\Big|\\
&&\leq \varepsilon\int_Q  e^{2\lambda\nu}\lambda^3\mu^4\beta^3\psi^2 dxdt
+C \int_0^T \int_{\tilde\omega_0} |f_0|^2_{C(\overline{Q})} e^{2\lambda\nu}\lambda^3\mu^4\beta^3\varphi^2dxdt;
\end{eqnarray*}
\begin{eqnarray*}
&&|I_2|=\Big|\int_Q e^{2\lambda\nu}\lambda^3\mu^4\beta^3\eta_1\varphi\psi_t dxdt\Big|\\
&&\leq \varepsilon \int_Q e^{2\lambda\nu}\frac{1}{\lambda\beta}\psi_t^2dxdt+
C\int^T_0\int_{\tilde\omega_0} e^{2\lambda\nu}\lambda^7\mu^8\beta^7\varphi^2 dxdt;
\end{eqnarray*}
\begin{eqnarray*}
&&|I_3|=\Big|\int_Q  \varphi\psi\eta_1(e^{2\lambda\nu}\lambda^3\mu^4\beta^3)_t dxdt\Big|
\leq C\Big|\int_Q e^{2\lambda\nu}\lambda^5\mu^4\beta^5\eta_1\varphi\psi dxdt\Big| \\
&&\leq \varepsilon\int_Q e^{2\lambda\nu}\lambda^3\mu^4\beta^3\psi^2 dxdt
+C\int^T_0\int_{\tilde\omega_0} e^{2\lambda\nu}\lambda^7\mu^4\beta^7\varphi^2 dxdt;
\end{eqnarray*}
\begin{eqnarray*}
&&|I_4|=\Big|\int_Q \sum\limits^n_{j=1} f_j\psi_{x_j} e^{2\lambda\nu}\lambda^3\mu^4
\beta^3\eta_1\varphi dxdt\Big|\leq C\sum\limits^n_{j=1} |f_j|_{C(\overline{Q})}
\int_Q e^{2\lambda\nu}\lambda^3\mu^4\beta^3\eta_1|\varphi||\nabla\psi|dxdt\\
&&\leq \varepsilon\int_Q e^{2\lambda\nu}\lambda\mu^2\beta|\nabla \psi|^2 dxdt
+C\int^T_0\int_{\tilde\omega_0}\sum\limits^n_{j=1} |f_j|^2_{C(\overline{Q})}
e^{2\lambda\nu}\lambda^5\mu^6\beta^5\varphi^2 dxdt;
\end{eqnarray*}
\begin{eqnarray*}
&&|I_5|=\Big|\int_Q\sum\limits^n_{j=1} f_j(e^{2\lambda\nu}\lambda^3\mu^4\beta^3\eta_1)_{x_j} \varphi\psi dxdt\Big|
\leq C\int^T_0\int_{\tilde\omega_0} \sum\limits^n_{j=1}|f_j|_{C(\overline{Q})}e^{2\lambda\nu}
\lambda^4\mu^5\beta^4|\varphi\psi|dxdt\\
&&\leq \varepsilon\int_Q e^{2\lambda\nu}\lambda^3\mu^4\beta^3\psi^2 dxdt
+C\sum\limits^n_{j=1}|f_j|^2_{C(\overline{Q})}\int^T_0\int_{\tilde\omega_0} e^{2\lambda\nu}
\lambda^5\mu^6\beta^5\varphi^2 dxdt;
\end{eqnarray*}

\begin{eqnarray*}
&&|I_6|=\Big|
\sum\limits^n_{i,j=1}\int_Q b^{ij}\varphi\psi(e^{2\lambda\nu}\lambda^3\mu^4\beta^3\eta_1)_{x_ix_j}dxdt\Big|\leq C\sum\limits^n_{i,j=1}|b^{i j}|_{C(\overline{Q})}\int_0^T\int_{\tilde\omega_0} 
e^{2\lambda\nu}\lambda^5\mu^6\beta^5|\varphi\psi|dxdt\\
&&\leq\varepsilon\int_Q e^{2\lambda\nu}\lambda^3\mu^4\beta^3\psi^2 dxdt
+C\sum\limits^n_{i,j=1}|b^{i j}|^2_{C(\overline{Q})}\int_0^T\int_{\tilde\omega_0} 
e^{2\lambda\nu}\lambda^7\mu^8\beta^7\varphi^2 dxdt;
\end{eqnarray*}

\begin{eqnarray*}
&&|I_7|=\Big|\sum\limits^n_{i,j=1}\int_Q (b^{ij}\psi)_{x_i}
(e^{2\lambda\nu}\lambda^3\mu^4\beta^3\eta_1)_{x_j}\varphi dxdt\Big|\\
&&\leq C\sum\limits^n_{i,j=1}|b^{i j}|_{C^{1, 0}(\overline{Q})} \int^T_0\int_{\tilde\omega_0}
e^{2\lambda\nu}\lambda^4\mu^5\beta^4|\varphi|(|\psi|+|\nabla\psi|)dxdt\\
&&\leq \varepsilon\int_Q e^{2\lambda\nu}\lambda\mu^2\beta|\nabla\psi|^2dxdt
+\varepsilon\int_Q e^{2\lambda\nu}\lambda^3\mu^4\beta^3 \psi^2dxdt\\
&&\quad+
C\int^T_0\int_{\tilde\omega_0} \sum\limits^n_{i,j=1}|b^{i j}|^2_{C^{1, 0}(\overline{Q})}
e^{2\lambda\nu}\lambda^7\mu^8\beta^7\varphi^2 dxdt;
\end{eqnarray*}

\begin{eqnarray*}
&&|I_8|=\Big|\sum\limits^n_{i,j=1}\int_Q  (b^{ij} \psi_{x_i})_{x_j} \varphi e^{2\lambda\nu}\lambda^3\mu^4\beta^3\eta_1 dxdt\Big|\\
&&\leq\varepsilon\int_Q e^{2\lambda\nu}\frac{1}{\lambda\beta}
\sum\limits^n_{i,j=1}  |\psi_{x_i x_j}|^2dxdt
+\varepsilon\int_Q e^{2\lambda\nu}\lambda\mu^2\beta
|\nabla\psi|^2dxdt\\
&&\quad\quad
+C\int^T_0\int_{\tilde\omega_0} \sum\limits^n_{i,j=1} |b^{i j}|^2_{C^{1, 0}(\overline{Q})} e^{2\lambda\nu}\lambda^7\mu^8\beta^7\varphi^2 dxdt;\\
\end{eqnarray*}
and
\begin{eqnarray*}
&&|I_9|=\Big|\sum\limits^n_{i,j=1}\int_Q b^{ij}\psi_{x_j}\varphi(e^{2\lambda\nu}\lambda^3\mu^4\beta^3\eta_1)_{x_i} dxdt\Big|\!\leq\! C\!\sum\limits^n_{i,j=1} |b^{i j}|_{C(\overline{Q})} \int_0^T\!\int_{\tilde\omega_0} \!
e^{2\lambda\nu}\lambda^4\mu^5\beta^4|\nabla\psi||\varphi|dxdt\\
&&\leq \varepsilon\int_Q e^{2\lambda\nu}\lambda\mu^2\beta|\nabla\psi|^2dxdt
+C\int^T_0\int_{\tilde\omega_0} \sum\limits^n_{i,j=1} |b^{i j}|^2_{C(\overline{Q})}
e^{2\lambda\nu}\lambda^7\mu^8\beta^7\varphi^2 dxdt.
\end{eqnarray*}

By the estimates from $I_1$ to $I_9$,   we get that for any $\varepsilon>0$, 
\begin{eqnarray*}
&&\int_Q e^{2\lambda\nu}\lambda^3\mu^4 \beta^3\eta_1\psi^2  dxdt\\
&&\leq \varepsilon\int_Q  e^{2\lambda\nu}(\lambda^3\mu^4\beta^3\psi^2+\lambda\mu^2\beta|\nabla\psi|^2)dxdt 
+\varepsilon\int_Q e^{2\lambda\nu}\frac{1}{\lambda\beta} \Big[\psi_t^2  +\sum\limits^n_{i,j=1}|\psi_{x_ix_j}|^2 \Big]  dxdt\\
&&\quad+C\Big(1+|f_0|^2_{C(\overline{Q})}+
\sum\limits_{j=1}^n|f_j|^2_{C(\overline{Q})}+\sum\limits_{i, j=1}^n|b^{i j}|^2_{C^{1, 0}(\overline{Q})}
\Big)\int_0^T \int_{\tilde\omega_0}e^{2\lambda\nu}\lambda^7\mu^8\beta^7\varphi^2   dxdt.
\end{eqnarray*}
Combining the above inequality with $(\ref{DX9})$ and choosing a sufficiently small $\varepsilon$, we have
\begin{eqnarray*}
\begin{array}{rl}
&\displaystyle\int_Q e^{2\lambda\nu}(\lambda\mu^2\beta|\nabla \varphi|^2+\lambda^3\mu^4\beta^3\varphi^2)dxdt
+\int_Q e^{2\lambda\nu}\frac{1}{\lambda\beta}\Big[\varphi_t^2
+\sum\limits^n_{i,j=1}|\varphi_{x_ix_j}|^2\Big]  
 dxdt \\
 &\displaystyle+\int_Q e^{2\lambda\nu}(\lambda\mu^2\beta|\nabla \psi|^2+\lambda^3\mu^4\beta^3\psi^2)dxdt
 +\int_Q e^{2\lambda\nu}\frac{1}{\lambda\beta}
 \Big[\psi_t^2+\sum\limits^n_{i,j=1}|\psi_{x_ix_j}|^2\Big] dxdt\\
 &\displaystyle\leq C(\mathcal{L})\int_0^T \int_{\tilde\omega_0}e^{2\lambda\nu}\lambda^7\mu^8\beta^7\varphi^2   dxdt.
\end{array}
\end{eqnarray*}
This completes the proof of (\ref{w}).
 \endpf

\medskip

\medskip

Finally,  we present an observability inequality for  the coupled system $(\ref{s})$. To this purpose,  
introduce the following functions:
\begin{eqnarray*}
&&l(t)=\begin{cases}
\ T^2/4    \qquad &  \mbox{ for } \text{$0\leq t \leq T/2$} ,\\[2mm]
\  t(T-t)      \qquad     &  \mbox{ for }\text{$T/2\leq t \leq T$},
\end {cases}\\
&&\bar{\nu}(x,t)=\frac{e^{\mu\eta(x)}-e^{2\mu|\eta|_{C(\overline{\Omega})}}}{l(t)}, 
 \quad \quad \bar{\beta}_0(t)=\frac{1}{l(t)},
\quad\quad \bar\beta(x, t)=\frac{e^{\mu\eta(x)}}{l(t)},
\end{eqnarray*}
\begin{equation}\label{notation}
\bar{\nu}_*(t)=\min\limits_{x \in \Omega}\bar{\nu}(x,t)\quad \mbox{ and }\quad  \hat{\rho}(t)=e^{-\lambda\bar{\nu}_*(t)}.
\end{equation}
\begin{proposition}\label{prop3}
Assume that  $\tilde\omega_0\cap \omega'\neq\emptyset$. 
Then there exist positive constants $C$ and $\tilde\mu_0$,  such that the following estimate  holds for 
any solution to $(\ref{s})$ and $\mu_i\geq\tilde\mu_0$ for $i=1, 2$:
$$
\int_\Omega\varphi^2(x, 0)  dx+ \int_Q \hat{\rho}^{-2}\theta_1^2 dxdt
+ \int_Q \hat{\rho}^{-2}\theta_2^2 dxdt
\leq C(\mathcal{L})\int_0^T \int_{\tilde\omega_0}e^{2\lambda\nu} \beta^7\varphi^2 dxdt,
$$
where   $\lambda$ and $\mu$ are fixed constants, and satisfy the conditions  in Proposition $\ref{prop2}$.
\end{proposition}

\noindent {\bf Proof. }  The whole proof is divided into three  parts. 

\medskip

{\bf Step 1. } First,  let $\eta_2 \in C^\infty([0,T])$ satisfy that
\begin{equation*}\label{b1}
0\leq \eta_2\leq 1 \mbox{ in }[0, T], \quad\eta_2 =1 \mbox{ in }  [0, T/2] \quad\mbox{ and }\quad  \eta_2 =0 \mbox{ in } [3T/4, T].
\end{equation*}
We multiply the first equation of $(\ref{t})$ by $\eta_2\varphi$ and integrate the associated equality in $\Omega$.
It follows that
\begin{eqnarray*}
&&\frac{1}{2}\frac{d}{dt}\int_{\Omega} \eta_2\varphi^2 dx\\
&&=\frac{1}{2}\int_{\Omega}\eta_{2, t}\varphi^2dx+
\int_{\Omega}\eta_2\varphi\Big[
-\sum\limits_{i, j=1}^{n}(b^{i j}\varphi_{x_i})_{x_j}-\sum\limits_{j=1}^{n}(f_j\varphi)_{x_j}+f_0\varphi-\xi_*\psi
\Big]dx\\
&&=\frac{1}{2}\int_{\Omega}\eta_{2, t}\varphi^2dx+\sum\limits_{i, j=1}^{n}\int_{\Omega}b^{i j}\varphi_{x_i}\varphi_{x_j}\eta_2 dx
+\int_{\Omega}\sum\limits_{j=1}^{n}\eta_2\varphi_{x_j} f_j\varphi dx
+\int_{\Omega}\eta_2 f_0\varphi^2dx-\int_{\Omega}\xi_* \eta_2\varphi\psi dx.
\end{eqnarray*}
This implies that
\begin{eqnarray*}
&&-\frac{d}{dt}\int_{\Omega} \eta_2\varphi^2 dx+\rho_0\int_{\Omega}|\nabla\varphi|^2\eta_2 dx\\
&&\leq -\int_{\Omega}\eta_{2, t}\varphi^2 dx+\int_{\Omega}\eta_2\psi^2 dx
+C\Big(1+|f_0|_{C(\overline{Q})}+
\sum\limits_{j=1}^{n} |f_j|^2_{C(\overline{Q})}\Big)\int_{\Omega}\eta_2\varphi^2 dx.
\end{eqnarray*}
For simplicity,  set
$$
b_0=C\Big(1+|f_0|_{C(\overline{Q})}+
\sum\limits_{j=1}^{n} |f_j|^2_{C(\overline{Q})}\Big).
$$
Multiplying the above inequality  by $e^{b_0t}$ and integrating the associated formula in $[0, T]$, 
  we get
\begin{eqnarray*}
&&\int_\Omega  \varphi^2(x, 0)  dx
+\rho_0 \int_0^{\frac{T}{2}}\int_\Omega   |\nabla\varphi|^2         dxdt\\
&&\leq \int_Q -e^{b_0 t}\eta_{2, t}\varphi^2 dxdt+e^{b_0 T}\int_Q \eta_2 \psi^2 dxdt\\
&&\leq Ce^{b_0 T}\int^{\frac{3T}{4}}_{\frac{T}{2}}\int_{\Omega} \varphi^2 dxdt+
Ce^{b_0 T}\int^{\frac{3T}{4}}_{0}\int_{\Omega} \psi^2 dxdt.
\end{eqnarray*}
By the Poincar\'e inequality, it follows that
\begin{eqnarray*}
&&\int_\Omega  \varphi^2(x, 0)  dx+
\int_0^{\frac{T}{2}}\int_\Omega (\varphi^2+|\nabla\varphi|^2)dxdt\\
&&\leq Ce^{b_0 T}\int^{\frac{3T}{4}}_{\frac{T}{2}} \int_{\Omega}(\varphi^2+\psi^2) dxdt+
Ce^{b_0 T}\int^{\frac{T}{2}}_{0}\int_{\Omega} \psi^2 dxdt.
\end{eqnarray*}

By the definitions of  $\bar{\nu}$ and $\bar{\beta}_0$ in $\Omega\times [0, T/2]$, 
$$
\displaystyle e^{2\lambda\bar\nu}=e^{2\lambda \frac{e^{\mu \eta}-e^{2\mu|\eta|_{C(\overline{\Omega})}}}{T^2/4} }
\leq  e^{8\lambda e^{\mu|\eta|_{C(\overline{\Omega})}}/T^2}\leq e^{8\lambda^2/T^2}
$$
and
$$
\displaystyle e^{2\lambda\bar\nu}=e^{2\lambda \frac{e^{\mu \eta}-e^{2\mu|\eta|_{C(\overline{\Omega})}}}{T^2/4} }
\geq e^{-8\lambda \frac{e^{2\mu|\eta|_{C(\overline{\Omega})}}}{T^2}}
\geq e^{-8\lambda^2/T^2}.
$$
Therefore,
\begin{eqnarray}\label{DX11}
  \begin{array}{rl}
&\displaystyle\int_\Omega  \varphi^2(x, 0) dx
+\int_0^{\frac{T}{2}}\int_\Omega   \lambda^3\mu^4\bar{\beta}_0^3e^{2\lambda\bar{\nu}}\varphi^2  dxdt
+\int_0^{\frac{T}{2}}\int_\Omega  \lambda\mu^2\bar{\beta}_0 e^{2\lambda\bar{\nu}}|\nabla\varphi|^2    dxdt\\[4mm]
&\displaystyle\leq  C \lambda^3\mu^4e^{b_0 T}e^{8\lambda^2/T^2}\Big[\int_{\frac{T}{2}}^{\frac{3}{4}T}\int_\Omega
  (\varphi^2  +  \psi^2) dxdt
+\int_0^{\frac{T}{2}}\int_\Omega \psi^2   dxdt\Big].
\end{array}
\end{eqnarray}

{\bf Step 2. }
In order to eliminate the terms related to $\psi$ in the right
 side of (\ref{DX11}), we multiply the second equation in $(\ref{t})$ by $\psi$ and  get that
\begin{eqnarray*}
&&\frac{1}{2}\frac{d}{dt} \int_\Omega \psi^2  dx
+\int_\Omega \sum\limits^n_{i, j=1} B^{i j}\psi_{x_i}\psi_{x_j}  dx
-\int_\Omega \sum\limits^n_{j=1} g_j\psi_{x_j}\psi dx-\int_{\Omega}g_0\psi^2 dx\\
&&\quad+\int_\Omega  (\frac{\nu_1}{\mu_1}\xi_{1}^2+\frac{\nu_2}{\mu_2}\xi^2_{2})\varphi\psi dx=0.
\end{eqnarray*}
It follows that
\begin{eqnarray*}
&&\frac{d}{dt}\int_{\Omega} \psi^2 dx+\rho_0\int_{\Omega}|\nabla \psi|^2 dx\\
&&\leq C\Big(1+\sum\limits_{j=1}^n |g_j|^2_{C(\overline{Q})}
+|g_0|_{C(\overline{Q})}\Big)\int_{\Omega} \psi^2 dx+\Big(\frac{1}{\mu_1^2}+\frac{1}{\mu_2^2}\Big)\int_{\Omega}\varphi^2 dx.
\end{eqnarray*}
For simplicity, set
$$
\tilde b_0=C\Big(1+|g_0|_{C(\overline{Q})}+
\sum\limits_{j=1}^{n} |g_j|^2_{C(\overline{Q})}\Big).
$$
By the Gronwall inequality, we get that for any $t\in [0, T/2]$,
$$
\displaystyle\int_\Omega \psi^2(x, t) dx
\leq  e^{\tilde b_0 T/2}
\Big(\frac{1}{\mu_1^2}+\frac{1}{\mu_2^2}\Big)
 \int_0^{\frac{T}{2}}\int_\Omega \varphi^2   dxdt.
$$

This,  together with (\ref{DX11}), implies that
\begin{eqnarray*}
  \begin{array}{rl}
&\displaystyle\int_\Omega  \varphi^2(x, 0) dx
+\int_0^{\frac{T}{2}}\int_\Omega   \lambda^3\mu^4\bar{\beta}_0^3e^{2\lambda\bar{\nu}}\varphi^2  dxdt
+\int_0^{\frac{T}{2}}\int_\Omega  \lambda\mu^2\bar{\beta}_0 e^{2\lambda\bar{\nu}}|\nabla\varphi|^2    dxdt\\[4mm]
&\displaystyle\leq  C \lambda^3\mu^4e^{b_0 T}e^{8\lambda^2/T^2}\Big[\int_{\frac{T}{2}}^{\frac{3}{4}T}\int_\Omega
  (\varphi^2  +  \psi^2) dxdt
+\frac{T}{2} e^{\tilde b_0 T/2}
\Big(\frac{1}{\mu_1^2}+\frac{1}{\mu_2^2}\Big)
 \int_0^{\frac{T}{2}}\int_\Omega \varphi^2   dxdt
\Big].
\end{array}
\end{eqnarray*}
Hence,  when 
\begin{equation}\label{mu}
\mu_k\geq \tilde\mu_0\deq Ce^{8\lambda^2/T^2+b_0 T/2+\tilde b_0 T/4}+\mu^*, \quad k=1, 2,
\end{equation}
where $\mu^*$ is the constant in (\ref{XX2}),  by (\ref{w}),  it holds that
\begin{eqnarray}\label{DX12}
  \begin{array}{rl}
&\displaystyle\int_\Omega  \varphi^2(x, 0) dx
+\int_0^{\frac{T}{2}}\int_\Omega   \lambda^3\mu^4\bar{\beta}_0^3e^{2\lambda\bar{\nu}}\varphi^2  dxdt
+\int_0^{\frac{T}{2}}\int_\Omega  \lambda\mu^2\bar{\beta}_0 e^{2\lambda\bar{\nu}}|\nabla\varphi|^2    dxdt\\[4mm]
&\displaystyle\leq  C \lambda^3\mu^4e^{b_0 T}e^{C\lambda^2}\int_{\frac{T}{2}}^{\frac{3}{4}T}\int_\Omega
  (\varphi^2  +  \psi^2) dxdt\\[4mm]
 &
 \displaystyle\leq C \lambda^3\mu^4e^{b_0 T}e^{C\lambda^2}\int_{\frac{T}{2}}^{\frac{3}{4}T}\int_\Omega
 e^{2\lambda\nu}\beta^3 (\varphi^2  +  \psi^2) dxdt\\[4mm]
 &\displaystyle\leq C(\mathcal{L}) \lambda^7\mu^8 e^{b_0 T}e^{C\lambda^2}\int^T_0\int_{\tilde\omega_0} e^{2\lambda\nu} \beta^7
\varphi^2dxdt.
\end{array}
\end{eqnarray}

{\bf \bf Step 3. }  Notice that 
$\bar{\nu}=\nu$ and $\bar{\beta}=\beta$ in $\Omega \times [T/2, T]$.
By the Carleman estimate (\ref{w}),   we obtain that
\begin{equation}\label{e1}
  \begin{split}
  &\int_{\frac{T}{2}}^T \int_\Omega   \lambda^3\mu^4\bar{\beta}^3_0 e^{2\lambda\bar{\nu}}\varphi^2   dxdt
+\int_{\frac{T}{2}}^T\int_\Omega  \lambda \mu^2\bar{\beta}_0 e^{2\lambda\bar{\nu}}|\nabla\varphi|^2       dxdt\\
&\leq\int_{\frac{T}{2}}^T \int_\Omega   \lambda^3\mu^4\bar{\beta}^3e^{2\lambda\bar{\nu}}\varphi^2   dxdt
+\int_{\frac{T}{2}}^T\int_\Omega  \lambda \mu^2\bar{\beta}e^{2\lambda\bar{\nu}}|\nabla\varphi|^2       dxdt\\
&\leq  C(\mathcal{L}) \int_0^T\int_{\tilde\omega_0}e^{2\lambda\nu}\lambda^7\mu^8\beta^7 \varphi^2   dxdt.
\end{split}
\end{equation}
Combining $(\ref{DX12})$ and $(\ref{e1})$, we get that
\begin{equation}\label{f1}
  \begin{split}
&\displaystyle\int_\Omega  \varphi^2(x, 0) dx
+\int_0^{T}\int_\Omega   \lambda^3\mu^4\bar{\beta}_0^3e^{2\lambda\bar{\nu}}\varphi^2  dxdt
+\int_0^{T}\int_\Omega  \lambda\mu^2\bar{\beta}_0 e^{2\lambda\bar{\nu}}|\nabla\varphi|^2    dxdt\\[4mm]
&\displaystyle\leq C(\mathcal{L}) \lambda^7\mu^8 e^{b_0 T}e^{C\lambda^2}\int^T_0\int_{\tilde\omega_0} e^{2\lambda\nu} \beta^7
\varphi^2dxdt.
\end{split}
\end{equation}

Notice that $\hat{\rho}(\cdot)$ is a non-decreasing strictly positive function blowing up at $t=T$.
Multiplying the second equation of (\ref{s}) by $\hat{\rho}^{-2}\theta_1$, 
using  integration by parts, we get
\begin{equation}
  \begin{split}
  \nonumber
\frac{1}{2}\frac{d}{dt} \int_\Omega   \hat{\rho}^{-2}\theta_1^2  dx
-\frac{1}{2}\int_\Omega  (\hat{\rho}^{-2})_t \theta_1^2  dx
+\int_\Omega   \sum\limits^n_{i, j=1}B^{i j}\theta_{1, x_i}\theta_{1, x_j} \hat{\rho}^{-2} dx\\
=\int_\Omega    -\frac{1}{\mu_1}\xi_{1}^2\hat{\rho}^{-2}\theta_1\varphi  dx
+\int_\Omega  g_0\hat{\rho}^{-2}\theta_1^2dx 
+\int_\Omega   \sum\limits^n_{j=1}g_j\theta_{1, x_j}\theta_1\hat{\rho}^{-2} dx.
\end{split}
\end{equation}
Since $(\hat{\rho}^{-2})_t\leq 0$, it follows that
\begin{equation}
  \begin{split}
  \nonumber
&\frac{1}{2}\frac{d}{dt}\int_\Omega  \hat{\rho}^{-2}\theta_1^2  dx
 +\rho_0\int_\Omega   |\nabla\theta_1|^2\hat{\rho}^{-2} dx\\
&\leq \frac{1}{2\mu_1^2}\int_\Omega  \varphi^2\hat{\rho}^{-2}   dx
+\frac{1}{2}\int_\Omega \hat{\rho}^{-2} \theta_1^2    dx
+\int_\Omega  |g_0|_{C(\overline{Q})} \theta_1^2\hat{\rho}^{-2} dx \\
& \q+C\sum\limits_{j=1}^{n} \int_\Omega  |g_j|_{C(\overline{Q})}^2\hat{\rho}^{-2}  \theta_1^2  dx
+\frac{\rho_0}{2}\int_\Omega    \hat{\rho}^{-2} |\nabla\theta_1|^2 dx.
\end{split}
\end{equation}
Hence, 
\begin{eqnarray*}
\begin{array}{ll}
\displaystyle\frac{d}{dt}\int_\Omega  \hat{\rho}^{-2} \theta_1^2  dx
+\int_\Omega   |\nabla\theta_1|^2\hat{\rho}^{-2} dx&\\
\displaystyle\leq  \frac{C}{\mu_1^2}\int_\Omega  \hat{\rho}^{-2}  \varphi^2 dx
+C\Big(1\!+\!|g_0|_{C(\overline{Q})}\!+\!\sum\limits_{j=1}^{n}|g_j|_{C(\overline{Q})}^2\Big)
\int_\Omega  \hat{\rho}^{-2} \theta_1^2  dx=
\frac{C}{\mu_1^2}\int_\Omega  \hat{\rho}^{-2}  \varphi^2 dx
+\tilde b_0
\int_\Omega  \hat{\rho}^{-2} \theta_1^2  dx.&
\end{array}
\end{eqnarray*}
 By the Gronwall inequality, we arrive at
\begin{equation}\label{DX15}
  \begin{split}
\int_\Omega  
\hat{\rho}^{-2}(t) \theta_1^2(x, t) dx
\leq    \frac{C}{\mu^2_1} e^{\tilde b_0 T}\int_Q \hat\rho^{-2}\varphi^2dxdt,\quad\forall\ t\in [0, T].
\end{split}
\end{equation}
In the same way, we can get
\begin{equation}\label{DX16}
  \begin{split}
\int_\Omega  
\hat{\rho}^{-2}(t) \theta_2^2(x, t) dx
\leq    \frac{C}{\mu^2_2} 
e^{\tilde b_0 T}
\int_Q \hat\rho^{-2}\varphi^2dxdt,\quad\forall\ t\in [0, T].
\end{split}
\end{equation}

By  (\ref{DX15}), (\ref{DX16}) and  (\ref{f1}),  for $\lambda$ and $\mu$ satisfying the conditions in Proposition \ref{prop2}, 
\begin{eqnarray*}
&&\int_Q
\hat{\rho}^{-2} \theta_1^2 dxdt+\int_Q  
\hat{\rho}^{-2} \theta_2^2 dxdt\\
&&\leq C\Big(\frac{1}{\mu_1^2}+\frac{1}{\mu_2^2}\Big)
e^{\tilde b_0 T}
\int_Q \hat\rho^{-2}\varphi^2dxdt\leq Ce^{\tilde b_0 T}
\int_Q e^{2\lambda\bar\nu}\varphi^2dxdt
\\
 &&\leq C(\mathcal{L})e^{(\tilde b_0+b_0) T} \lambda^4\mu^4 e^{C\lambda^2}\int^T_0\int_{\tilde\omega_0} e^{2\lambda\nu} \beta^7
\varphi^2dxdt\leq C(\mathcal{L})\int^T_0\int_{\tilde\omega_0} e^{2\lambda\nu} \beta^7
\varphi^2dxdt,
\end{eqnarray*}
where the constants $C(\mathcal{L})$ in different places may be different.  
The constant $C$ in $Ce^{2C\lambda^2}$ for the last term is sufficiently large, such that
 $e^{2C\lambda^2}$ is larger than $C e^{(\tilde b_0+b_0) T} \lambda^4\mu^4 e^{C\lambda^2}\mathcal{L}^3$. 
This, together with (\ref{f1}),  completes the proof of Proposition \ref{prop3}.
\endpf

\section{Null controllability of coupled  linear parabolic  systems}

In this section, we study the  controllability of a coupled linear parabolic system in the framework of classical solutions. The key is 
to improve  regularity of control functions in a suitable H\"older space. The considered  controlled system
 is  as follows:

\begin{equation}\label{s*}
  \begin{cases}
  \begin{split}
&\bar y_t-\sum\limits^n_{i, j=1}(b^{i j}\bar{y}_{x_i})_{x_j}+
\sum\limits_{j=1}^{n} f_j \bar{y}_{x_j}+f_0\bar y=\xi_0 u+\frac{1}{\mu_1} \xi^2_1 p_1+\frac{1}{\mu_2} \xi^2_2 p_2  &\mbox{ in }Q,&\\
&p_{1, t}+\sum\limits^n_{i, j=1}(B^{i j}p_{1, x_i})_{x_j}+
\sum\limits_{j=1}^{n} g_j p_{1, x_j}+g_0 p_1=\nu_1 \xi_*(\bar y-y_{1, d})   &\mbox{ in }  Q,&\\
&p_{2, t}+\sum\limits^n_{i,j=1}(B^{ij}p_{2, x_i})_{x_j}+\sum\limits_{j=1}^{n} g_j p_{2, x_j}+g_0 p_2=
\nu_2 \xi_*(\bar y-y_{2, d})    &\mbox{ in }  Q,&\\
&\bar y=p_1= p_2=0   &\mbox{ on } \Sigma,&\\
&\bar y(x, 0)=y_0(x),\  p_1(x, T)=p_2(x, T)=0 &\mbox{ in } \Omega,&
  \end{split}
  \end{cases}
\end{equation}
where  $b^{i j}$, $B^{i j}$, $f_j$, $g_j$, $f_0$ and $g_0$ satisfy  the same conditions in (\ref{s}) for $i, j=1,\cdots, n$, and 
$\xi_0$, $\xi_1$, $\xi_2$ and $\xi_*$ satisfy the same conditions in (\ref{a}) and (\ref{d}).

\medskip

We have the following partial null controllability result for the coupled linear parabolic system (\ref{s*})  with a control function
$u$ in the H\"older space $C^{1+\alpha, \frac{1+\alpha}{2}}(\overline{Q})$.
\begin{proposition}\label{prop4}Assume that all conditions in Proposition $\ref{prop3}$ hold. Then 
for any $y_0 \in C^{3+\alpha}_0(\Omega)$, there is a control
$u\in C^{1+\alpha, \frac{1+\alpha}{2}}(\overline{Q})$,  such that the corresponding solution $(\bar{y}, p_1, p_2)$ to $(\ref{s*})$ satisfies
$$\bar{y}(x, T)=0 \qquad in \q \Omega. $$
Moreover,
$$|u|_{C^{1+\alpha, \frac{1+\alpha}{2}}(\overline{Q})}\leq C(\mathcal{L})\Big(|y_0|_{L^2(\Omega)}
+|\hat\rho y_{1, d}|_{L^2(Q)}+|\hat\rho y_{2, d}|_{L^2(Q)}\Big).$$
\end{proposition}

\noindent {\bf Proof. } The whole proof is divided into three parts. The ideas originate from \cite{bar} and \cite{Liux}.

\medskip

{\bf Step 1. }
First,  for any  $\varepsilon>0$,  construct the following functional:
$$J_{\varepsilon}(u)=\frac{1}{2}\int_{Q} e^{-2\lambda\nu}\beta^{-7}u^2dxdt
+\frac{1}{2\varepsilon}\int_{\Omega} \bar{y}^2(x, T) dx, \quad\forall\ u\in\widetilde{\mathcal{U}},$$
where $\displaystyle\widetilde{\mathcal{U}}=
\Big\{\ u\in L^2(Q)\  \Big|\  \int_{Q}e^{-2\lambda\nu}\beta^{-7}u^2dxdt<\infty\  \Big\}$ 
and $\bar y\in C([0, T]; L^2(\Omega))$ is the component of  solution to (\ref{s*}) corresponding to 
$u\in\mathcal{U}$.

By the variational method,    $J_{\varepsilon}(\cdot)$ has
a unique minimum element $u_{\varepsilon}\in L^2(Q)$. Denote by 
$(\bar y_{\varepsilon}, p_{1}^\varepsilon, p_2^\varepsilon)$ the  solution to 
(\ref{s*}) associated to $u=u_\varepsilon$.  Then,  it follows that
\begin{equation}\label{o1}
u_{\varepsilon}=e^{2\lambda\nu}\beta^{7}\xi_0\varphi_\varepsilon,
\end{equation}
where $\varphi_\varepsilon$ is the component of  solution
 $(\varphi_\varepsilon, \theta_1^\varepsilon, \theta_2^\varepsilon)$ to the following coupled system:
\begin{equation}\label{s**}
  \begin{cases}
  \begin{split}
&\varphi_{\varepsilon, t}+\sum\limits^n_{i,j=1}(b^{ij}\varphi_{\varepsilon, x_i})_{x_j}+
\sum\limits_{j=1}^{n} (f_j\varphi_\varepsilon)_{x_j}-f_0\varphi_\varepsilon+
 \xi_{*}(\nu_1\theta_1^\varepsilon+\nu_2\theta_2^\varepsilon)=0  &\mbox{ in }Q,&\\
&\theta_{1, t}^\varepsilon-\sum\limits^n_{i,j=1}(B^{ij}\theta^\varepsilon_{1, x_i})_{x_j}+
\sum\limits_{j=1}^{n} (g_j\theta_1^\varepsilon)_{x_j}-g_0\theta_1^\varepsilon+
 \frac{1}{\mu_1}\xi_{1}^2\varphi_\varepsilon=0   &\mbox{ in }  Q,&\\
&\theta^\varepsilon_{2, t}-\sum\limits^n_{i,j=1}(B^{ij}\theta^\varepsilon_{2, x_i})_{x_j}+
\sum\limits_{j=1}^{n} (g_j\theta^\varepsilon_2)_{x_j}-g_0\theta_2^\varepsilon+
 \frac{1}{\mu_2}\xi_{2}^2\varphi_\varepsilon=0    &\mbox{ in }  Q,&\\
&\varphi_\varepsilon=\theta_1^\varepsilon= \theta_2^\varepsilon=0   &\mbox{ on } \Sigma,&\\
&\varphi_\varepsilon(x, T)=-\displaystyle\frac{1}{\varepsilon} \bar y_\varepsilon(x, T), \  \theta_1^\varepsilon(x, 0)=\theta_2^\varepsilon(x, 0)=0 &\mbox{ in } \Omega.&
  \end{split}
  \end{cases}
\end{equation}

Next, by the duality between $(\bar y_\varepsilon, p_1^\varepsilon, p_2^\varepsilon)$
and $(\varphi_\varepsilon, \theta^\varepsilon_1, \theta^\varepsilon_2)$
in (\ref{s*}) and (\ref{s**}),  according to Proposition \ref{prop3}, we  obtain that 
\begin{eqnarray*}
&&\frac{1}{\varepsilon} \int_{\Omega} \bar y_{\varepsilon}^2(x, T)dx
+\int_Q \xi_0\varphi_\varepsilon u_{\varepsilon}dxdt\\
&&=\frac{1}{\varepsilon} \int_{\Omega} \bar y_{\varepsilon}^2(x, T)dx
+\int_Q \xi_0^2 e^{2\lambda\nu}\beta^7 \varphi_\varepsilon^2 dxdt\\
&&=-\int_{\Omega} \varphi_\varepsilon(x, 0)y_0(x)dx
+\int_Q\xi_*(\nu_1\theta_1^\varepsilon y_{1, d}+\nu_2\theta_2^\varepsilon y_{2, d})
dxdt\\
&&\leq |\varphi_\varepsilon(\cdot, 0)|_{L^2(\Omega)}|y_0|_{L^2(\Omega)}
+C|\hat\rho^{-1} \theta_1^\varepsilon|_{L^2(Q)}|\hat\rho y_{1, d}|_{L^2(Q)}
+C|\hat\rho^{-1} \theta_2^\varepsilon|_{L^2(Q)}|\hat\rho y_{2, d}|_{L^2(Q)}\\[3mm]
&&\leq C\Big(|\varphi_\varepsilon(\cdot, 0)|_{L^2(\Omega)}+
|\hat\rho^{-1} \theta_1^\varepsilon|_{L^2(Q)}
+|\hat\rho^{-1} \theta_2^\varepsilon|_{L^2(Q)}\Big)\Big(|y_0|_{L^2(\Omega)}
+|\hat\rho y_{1, d}|_{L^2(Q)}+|\hat\rho y_{2, d}|_{L^2(Q)}\Big)\\
&&\leq C(\mathcal{L})\Big(\int^T_0\int_{\tilde\omega_0} e^{2\lambda\nu}
\beta^7 \varphi_\varepsilon^2 dxdt\Big)^{1/2}\Big(|y_0|_{L^2(\Omega)}
+|\hat\rho y_{1, d}|_{L^2(Q)}+|\hat\rho y_{2, d}|_{L^2(Q)}\Big).
\end{eqnarray*}
This implies that
\begin{eqnarray}\label{s1}
\begin{array}{ll}
&\displaystyle\frac{1}{\varepsilon} \int_{\Omega} \bar y_{\varepsilon}^2(x, T)dx
+\int_Q \xi_0^2 e^{2\lambda\nu}\beta^7 \varphi_\varepsilon^2 dxdt\displaystyle
\leq C(\mathcal{L})\Big(|y_0|_{L^2(\Omega)}
+|\hat\rho y_{1, d}|_{L^2(Q)}+|\hat\rho y_{2, d}|_{L^2(Q)}\Big)^2.
\end{array}
\end{eqnarray}

 \medskip
 
 {\bf Step 2. } In this step, we  prove that $u_\varepsilon$ in (\ref{o1}) 
 are  bounded in $C^{1+\alpha, \frac{1+\alpha}{2}}(\overline{Q})$ 
 with respect to $\varepsilon$.
 Put $\psi^\varepsilon=\nu_1\theta^\varepsilon_1+\nu_2\theta^\varepsilon_2$ and
 $(\ref{s**})$ is simplified  as
\begin{equation}\label{t1}
  \begin{cases}
  \begin{split}
&\varphi_{\varepsilon, t}+\sum\limits^n_{i,j=1}(b^{ij}\varphi_{\varepsilon, x_i})_{x_j}+
\sum\limits_{j=1}^{n} (f_j\varphi_\varepsilon)_{x_j}-f_0\varphi_\varepsilon+
 \xi_{*}\psi^\varepsilon=0  &\mbox{ in }Q,&\\
&\psi_{t}^\varepsilon-\sum\limits^n_{i,j=1}(B^{ij}\psi^\varepsilon_{x_i})_{x_j}+
\sum\limits_{j=1}^{n} (g_j\psi^\varepsilon)_{x_j}-g_0\psi^\varepsilon+
 \Big(\frac{\nu_1}{\mu_1}\xi_{1}^2+\frac{\nu_2}{\mu_2}\xi_{2}^2\Big)
 \varphi_\varepsilon=0   &\mbox{ in }  Q,&\\
&\varphi_\varepsilon=\psi^\varepsilon=0   &\mbox{ on } \Sigma,&\\
&\varphi_\varepsilon(x, T)=-\displaystyle\frac{1}{\varepsilon} \bar y_\varepsilon(x, T), 
\  \psi^\varepsilon(x, 0)=0 &\mbox{ in } \Omega.&
  \end{split}
  \end{cases}
\end{equation}

Let $\hat\delta_0$ be a positive number and $\{\delta_k\}_{k\geq1}$ be a monotone  increasing sequence
such that $0<\delta_k<\hat\delta_0<\lambda/2$. 
Recall that
$$
\beta_0(t)=\frac{1}{t(T-t)}\quad\mbox{ and }\quad
\nu_0(t)=(1-e^{2\mu|\eta|_{C(\overline{\Omega})}})\beta_0(t),\ 
$$
and
write
$$
m^k_\varepsilon=e^{(\lambda+\delta_k)\nu_0}\beta^7_0\varphi_\varepsilon \quad
\mbox{ and }\quad
n^k_\varepsilon=e^{(\lambda+\delta_k)\nu_0}\beta^7_0\psi^\varepsilon.
$$
Then it is easy to check that $m^k_\varepsilon$ and $n^k_\varepsilon$ satisfy
\begin{equation}\label{u1}
  \begin{cases}
  \begin{split}
&m^k_{\varepsilon,t}+
\sum\limits^n_{i,j=1}(b^{i j}m^k_{\varepsilon,x_i})_{x_j}+
\sum\limits_{j=1}^{n} (f_j m^k_\varepsilon)_{x_j}-f_0 m^k_\varepsilon\\
&\qq\qq\qq\qq\qq=-n^k_\varepsilon\xi_{*}+(e^{(\lambda+\delta_k)\nu_0}\beta^7_0)_t\varphi_\varepsilon\triangleq g^k_\varepsilon \qquad &\mbox{ in }Q,&\\
&n^k_{\varepsilon,t}-\sum\limits_{i, j=1}^n (B^{i j}n^k_{\varepsilon, x_i})_{x_j}
+
\sum\limits_{j=1}^{n} (g_j n^k_\varepsilon)_{x_j}
-g_0 n^k_\varepsilon\\
& \qq\qq\qq\qq\qq=
-\Big(\frac{\nu_1}{\mu_1}\xi_{1}^2+\frac{\nu_2}{\mu_2}\xi_{2}^2\Big)
m^k_\varepsilon+(e^{(\lambda+\delta_k)\nu_0}\beta^7_0)_t\psi^\varepsilon
 \triangleq l^k_\varepsilon\qquad &\mbox{ in }Q,&\\
&m^k_\varepsilon=n^k_\varepsilon=0 \qquad &\mbox{ on }\Sigma,&\\
&m^k_\varepsilon(x, 0)=m^k_\varepsilon(x, T)=0,\q n^k_\varepsilon(x, 0)
=n^k_\varepsilon(x, T)=0\qquad&\mbox{ in }\Omega.&
  \end{split}
  \end{cases}
\end{equation}
For $k=1$, 
$$ 
\ds  g^1_\varepsilon=-n^1_\varepsilon\xi_{*}+(e^{(\lambda+\delta_1)\nu_0}\beta^7_0)_t\varphi_\varepsilon
=-e^{(\lambda+\delta_1)\nu_0}\beta^7_0\psi^\varepsilon\xi_*+(e^{(\lambda+\delta_1)\nu_0}\beta^7_0)_t\varphi_\varepsilon.
$$
By the Carleman estimate (\ref{w})  for (\ref{t1}) and  (\ref{s1}),
noting that $e^{2\delta_1 \nu_0}\beta_0^{15}\leq C$, we get
\begin{eqnarray*}
&&|g^1_\varepsilon|^2_{L^2(Q)} \leq
\int_Q e^{2(\lambda+\delta_1)\nu_0}\beta_0^{14}|\psi^\varepsilon|^2dxdt+
\int_Q C\lambda^4 e^{2(\lambda+\delta_1)\nu_0}\beta_0^{18} \varphi_\varepsilon^2dxdt\\
&&\leq C\int_Q\lambda^4 e^{2\lambda\nu}\beta^3(|\psi^\varepsilon|^2+\varphi_\varepsilon^2)dxdt
\leq C(\mathcal{L})\int^T_0\int_{\tilde\omega_0} e^{2\lambda\nu}\lambda^8\mu^4\beta^7 \varphi_\varepsilon^2 dxdt\\
&&\leq C(\mathcal{L})\Big(|y_0|_{L^2(\Omega)}
+|\hat\rho y_{1, d}|_{L^2(Q)}+|\hat\rho y_{2, d}|_{L^2(Q)}\Big)^2,
\end{eqnarray*}
since $\lambda$ and $\mu$ depend on $\mathcal{L}$. Further, 
$$ l^1_\varepsilon
=-(\frac{\nu_1}{\mu_1}\xi^2_1+\frac{\nu_2}{\mu_2}\xi^2_2)e^{(\lambda+\delta_1)\nu_0}\beta^7_0\varphi_\varepsilon+
(e^{(\lambda+\delta_1)\nu_0}\beta^7_0)_t\psi^\varepsilon.
 $$
 Similarly,  by the Carleman estimate (\ref{w}) for (\ref{t1}) and  (\ref{s1}), we have
\begin{eqnarray*}
&&|l^1_\varepsilon|^2_{L^2(Q)}  
 \leq \int_Q Ce^{2(\lambda+\delta_1)\nu_0}\beta^{14}_0 \varphi_\varepsilon^2dxdt
+\int_Q Ce^{2(\lambda+\delta_1)\nu_0}\lambda^4\beta^{18}_0 |\psi^\varepsilon|^2dxdt\\
&&\leq C\int_Q\lambda^4 e^{2\lambda\nu}\beta^3(|\psi^\varepsilon|^2+\varphi_\varepsilon^2)dxdt
\leq C(\mathcal{L})\int^T_0\int_{\tilde\omega_0} e^{2\lambda\nu}\lambda^8\mu^4\beta^7 \varphi_\varepsilon^2 dxdt\\
&&\leq C(\mathcal{L})\Big(|y_0|_{L^2(\Omega)}
+|\hat\rho y_{1, d}|_{L^2(Q)}+|\hat\rho y_{2, d}|_{L^2(Q)}\Big)^2.
\end{eqnarray*}
By $L^p$-estimates for linear  parabolic equations of second order  (see  \cite[Lemma 4.1]{Liu}), we get that
\begin{eqnarray*}
&&|m^1_\varepsilon|^2_{W_2^{2,1}(Q)}
+|n^1_\varepsilon|^2_{W_2^{2,1}(Q)}
 \leq C(\mathcal{L})(|g^1_\varepsilon|^2_{L^2(Q)}+|l^1_\varepsilon|^2_{L^2(Q)})\\[2mm]
 &&\leq C(\mathcal{L})\Big(|y_0|_{L^2(\Omega)}
+|\hat\rho y_{1, d}|_{L^2(Q)}+|\hat\rho y_{2, d}|_{L^2(Q)}\Big)^2.
\end{eqnarray*}
By the Sobolev embedding, $W_2^{2,1}(Q)\hookrightarrow L^{r_1}(Q)$, for
\begin{equation}
\nonumber
r_1=\begin{cases}
2(n+2)/(n-2),  \qquad &  \mbox{if} \quad\text{$n>2$} ,\\
\mbox{any constant larger than }1,      \qquad     & \mbox{if} \quad\text{$n\leq2$}.
\end {cases}
\end{equation}
Hence, it follows that
\begin{eqnarray*}
&&\ds |m^1_\varepsilon|^2_{L^{r_1}(Q)} + |n^1_\varepsilon|^2_{L^{r_1}(Q) } 
 \leq C(|m^1_\varepsilon|^2_{W_2^{2,1}(Q)}+|n^1_\varepsilon|^2_{W_2^{2,1}(Q)})\\
&&\leq C(\mathcal{L})\Big(|y_0|_{L^2(\Omega)}
+|\hat\rho y_{1, d}|_{L^2(Q)}+|\hat\rho y_{2, d}|_{L^2(Q)}\Big)^2.
\end{eqnarray*}

{\bf Step 3. } In the following, we give the estimates for $m^2_\varepsilon$ and $n^2_\varepsilon$. Notice that
\begin{eqnarray*}
&&\ds  g^2_\varepsilon=-n^2_\varepsilon\xi_{*}+(e^{(\lambda+\delta_2)\nu_0}\beta^7_0)_t\varphi_\varepsilon
=-e^{(\lambda+\delta_2)\nu_0}\beta^7_0\psi^\varepsilon\xi_*+
(e^{(\lambda+\delta_2)\nu_0}\beta^7_0)_t\varphi_\varepsilon\\
&&=-e^{(\delta_2-\delta_1)\nu_0}n^1_\varepsilon\xi_*+
(e^{(\lambda+\delta_2)\nu_0}\beta^7_0)_t    e^{-(\lambda+\delta_1)\nu_0}\beta_0^{-7} m^1_\varepsilon,
\end{eqnarray*}
and
\begin{eqnarray*}
&&l^2_\varepsilon
=-(\frac{\nu_1}{\mu_1}\xi^2_1+\frac{\nu_2}{\mu_2}\xi^2_2)e^{(\lambda+\delta_2)\nu_0}\beta^7_0\varphi_\varepsilon+
(e^{(\lambda+\delta_2)\nu_0}\beta^7_0)_t\psi^\varepsilon\\
&&=-(\frac{\nu_1}{\mu_1}\xi^2_1+\frac{\nu_2}{\mu_2}\xi^2_2)e^{(\delta_2-\delta_1)\nu_0} m^1_\varepsilon+
(e^{(\lambda+\delta_2)\nu_0}\beta^7_0)_t e^{-(\lambda+\delta_1)\nu_0} \beta_0^{-7} n^1_\varepsilon.
\end{eqnarray*}
This implies that 
$
\displaystyle|g_\varepsilon^2|+|l^2_\varepsilon|\leq 
C\lambda^2(|m_\varepsilon^1|+|n_\varepsilon^1|).
$
Hence,
\begin{eqnarray*}
&&|g_\varepsilon^2|^2_{L^{r_1}(Q)}+|l_\varepsilon^2|^2_{L^{r_1}(Q)}
\leq C\lambda^4(|m_\varepsilon^1|^2_{L^{r_1}(Q)}+|n_\varepsilon^1|^2_{L^{r_1}(Q)})\\
&&\leq C(\mathcal{L})\Big(|y_0|_{L^2(\Omega)}
+|\hat\rho y_{1, d}|_{L^2(Q)}+|\hat\rho y_{2, d}|_{L^2(Q)}\Big)^2.
\end{eqnarray*}

Again, by $L^p$-estimates for linear  parabolic equations of second order,  it holds that
 $m^2_\varepsilon$, $n^2_\varepsilon \in W_{r_1}^{2,1}(Q)$. Moreover,
\begin{eqnarray*}
&&|m^2_\varepsilon|^2_{W_{r_1}^{2,1}(Q)}+|n^2_\varepsilon|^2_{W_{r_1}^{2,1}(Q)} 
\leq C(\mathcal{L}) (|g_\varepsilon^2|^2_{L^{r_1}(Q)}+|l_\varepsilon^2|^2_{L^{r_1}(Q)})\\
&&\leq C(\mathcal{L})\Big(|y_0|_{L^2(\Omega)}
+|\hat\rho y_{1, d}|_{L^2(Q)}+|\hat\rho y_{2, d}|_{L^2(Q)}\Big)^2.
\end{eqnarray*}
By the Sobolev embedding, we see that $W_{r_1}^{2,1}(Q)\hookrightarrow L^{r_2}(Q)$, where
\begin{equation}
\nonumber
r_2=\begin{cases}
r_1(n+2)/(n+2-2r_1), \qquad & \mbox{if} \quad\text{$n+2-2r_1>0$} ,\\
\mbox{any constant  larger than } 1,      \qquad     & \mbox{if} \quad\text{$n+2-2r_1\leq0$}.
\end {cases}
\end{equation}

We repeat the above arguments and define $r_N$  similarly, for any $N\in\mathbb N$. Note that $\{r_N\}_{N\in\mathbb N}$ 
is monotone increasing and
$
 \displaystyle\frac{1}{r_N}-\frac{1}{r_{N+1}}=\frac{2}{n+2}$, $\forall$ $N\in\mathbb{N}.$
Then there exists a positive integer $N^*$,  such that $r_{N^*}>(n+2)/(1-\alpha)$.
By the Sobolev embedding, $W_{r_{N^*}}^{2,1}(Q)\hookrightarrow 
C^{1+\alpha, \frac{1+\alpha}{2}}(\overline{Q})$.
Therefore,
\begin{eqnarray*}
\begin{split}
&|u_\varepsilon|^2_{C^{1+\alpha, \frac{1+\alpha}{2}}(\overline{Q})}
=\Big|e^{2\lambda\nu}\beta^{7} 
e^{-(\lambda+\delta_{N^*})\nu_0}\beta_0^{-7} m_\varepsilon^{N^*}
\xi_0\Big|_{C^{1+\alpha, \frac{1+\alpha}{2}}(\overline{Q})}^2\\
&  \leq C|m_\varepsilon^{N^*}|^2_{C^{1+\alpha, \frac{1+\alpha}{2}}(\overline{Q})}
\leq C|m_\varepsilon^{N^*}|^2_{W_{r_{N^*}}^{2,1}(Q)} \leq C(\mathcal{L})\Big(|y_0|_{L^2(\Omega)}
+|\hat\rho y_{1, d}|_{L^2(Q)}+|\hat\rho y_{2, d}|_{L^2(Q)}\Big)^2.
\end{split}
\end{eqnarray*}

	Finally, letting $\varepsilon\rightarrow 0$ in $(\ref{s1})$ and the above equality, we can find  a control function
$u\in C^{1+\alpha, \frac{1+\alpha}{2}}(\overline{Q})$,  such that the corresponding solution 
$\bar{y}$ to $(\ref{s*})$ satisfies
$\bar{y}(x, T)=0  \mbox{ in   }\Omega.$
Moreover,
$$|u|_{C^{1+\alpha, \frac{1+\alpha}{2}}(\overline{Q})}\leq\! C(\mathcal{L})\big(|y_0|_{L^2(\Omega)}
\!\!+\!|\hat\rho y_{1, d}|_{L^2(Q)}\!\!+\!|\hat\rho y_{2, d}|_{L^2(Q)}\big)^2.$$
This completes the proof. 
\endpf

\section{Proofs of main controllability results}
\subsection{Controllability of the coupled quasi-linear parabolic system $(\ref{dl2})$}

In this subsection, we prove a local  controllability result of the coupled 
 quasi-linear parabolic system $(\ref{dl2})$ by 
 the fixed point technique.  By Proposition \ref{DX!}, this controllability result implies that Theorem \ref{th1} 
 holds true.
 
 For this purpose,   for $y_0\in C^{3+\alpha}_0(\Omega)$,  define
 $$
 K=\Big\{\ z\in C^{3+\alpha, \frac{3+\alpha}{2}}(\overline{Q})\ \Big|\ 
 |z|_{C^{3+\alpha, \frac{3+\alpha}{2}}(\overline{Q})}\leq 1, z(x, 0)=y_0(x)\ \mbox{in }\Omega\ \Big\}.
 $$
Then $K$ is a nonempty convex and compact subset of $L^2(Q)$ for small initial datum $y_0\in C^{3+\alpha}_0(\Omega)$.
For any $z\in K$,   consider the following linearized parabolic equation of $(\ref{dl2})$:
\begin{eqnarray}\label{w1}
\left\{
\begin{array}{ll}
\bar y_{t}-\sum\limits^n_{i, j=1}\big(a^{ij}(z,\nabla z)\bar y_{ x_i}\big)_{x_j}&\\
\quad\quad\quad\quad\quad
+F_1(z, \nabla z)\bar y+F_2(z, \nabla z)\cdot \nabla \bar y=\xi_{0}u+
\displaystyle\frac{1}{\mu_1}
\xi_{1}^2 p_1
+\frac{1}{\mu_2}
\xi_{2}^2 p_2     &\mbox{ in }   Q,\\
p_{1, t}+\sum\limits^n_{i, j=1} \big(A^{ij}(z, \nabla z)
p_{1, x_i}\big)_{x_j} &\\
\quad\quad\quad\quad\quad+\sum\limits_{j=1}^n 
e^j(z, \nabla  z)  p_{1, x_j}+d_0(z, \nabla z)
p_1=
\nu_1 \xi_{*}(\bar y-y_{1, d}) &\mbox{ in }Q,\\
p_{2, t}+\sum\limits^n_{i, j=1} \Big(A^{ij}(z, \nabla z)
p_{2, x_i}\Big)_{x_j}&\\
\quad\quad\quad\quad\quad+\sum\limits_{j=1}^n 
e^j(z, \nabla z)  p_{2, x_j}+d_0(z, \nabla z)
p_2=
\nu_2 \xi_{*}(\bar y-y_{2, d}) &\mbox{ in }Q,\\
\bar y=p_1=p_2=0  &\mbox{ on }  \Sigma,\\[2mm]
\bar y(x, 0)=y_0(x), \ p_1(x,  T)=p_2(x, T)=0 &\mbox{ in }  \Omega,
\end{array}
\right.
\end{eqnarray}
where 
$$\displaystyle F_1(z, \nabla z)=\int^1_0\!\! f_y(sz, s\nabla z)ds\mbox{ and }
F_2(z, \nabla z)=(F_2^1(z, \nabla z), \cdots, F_2^n(z, \nabla z))^\top=\int^1_0\!\! \nabla_\zeta f(sz, s\nabla z)ds.$$

\medskip

First,  we have the following well-posedness  results for linear parabolic equations of second order (e. g. 
\cite[Theorem 5.2 on Page 320]{La} and \cite[Theorem 4.28 on Page 77]{lie}).
\begin{lemma}\label{lemma*}
Consider the following linear parabolic equation:
\begin{equation}\label{u**}
  \begin{cases}
  \begin{split}
  & v_t-\sum\limits_{i, j=1}^n b^{i j} v_{x_ix_j}+\sum\limits_{j=1}^n f_j v_{x_j}+f_0 v=h_0 \qquad &\mbox{ in }  Q,& \\
  & v=0 \qquad &\mbox{ on } \Sigma,&\\
  & v(x, 0)=v_0(x) \qquad  &\mbox{ in } \Omega,& \\
  \end{split}
  \end{cases}
\end{equation}
where $b^{i j}$ satisfy the conditions in $(\ref{s})$,  and 
$b^{i j}, f_j, f_0\in C^{1+\alpha, \frac{1+\alpha}{2}}(\overline{Q})$ for $i, j=1, \cdots, n$. 
 Then for any   $v_0\in C_0^{3+\alpha}(\Omega)$ and 
 $h_0\in C^{1+\alpha, \frac{1+\alpha}{2}}(\overline{Q})$ 
 with $h_0(x, 0)=0$ on $\partial\Omega$,  the equation $(\ref{u**})$ admits a unique solution
 $v\in C^{3+\alpha, \frac{3+\alpha}{2}}(\overline{Q})$. Moreover, 
 $$
 |v|_{C^{3+\alpha, \frac{3+\alpha}{2}}(\overline{Q})}
 \leq C(n, \Omega, T, \alpha, \rho_0, \Lambda)
 \Big(|h_0|_{C^{1+\alpha, \frac{1+\alpha}{2}}(\overline{Q})}
 +|v_0|_{C^{3+\alpha}(\overline{\Omega})}\Big),
 $$ 
 with $$\sum\limits_{i, j=1}^n |b^{ i j}|_{C^{1+\alpha, \frac{1+\alpha}{2}}(\overline{Q})}
 +\sum\limits_{j=1}^n |f_{j}|_{C^{1+\alpha, \frac{1+\alpha}{2}}(\overline{Q})}
 +|f_0|_{C^{1+\alpha, \frac{1+\alpha}{2}}(\overline{Q})}\leq \Lambda.$$
 \end{lemma}

By Proposition $\ref{prop4}$,  we have the following    controllability result for (\ref{w1}).
\begin{lemma}\label{lem4}
Assume that  $\tilde\omega_0\cap \omega'\neq \emptyset$. Then 
there exist positive constants $\widehat{L}$ and $\hat\mu_0$, such that for any
$y_0\in C^{3+\alpha}_0(\Omega)$ and $\mu_k\geq \hat\mu_0$ for $k=1, 2$, one can always find  a control
 $u\in C^{1+\alpha, \frac{1+\alpha}{2}}(\overline{Q})$  
so that the corresponding solution to the system $(\ref{w1})$ satisfies
$$\bar{y}(x, T)=0  \qquad \mbox{ in } \Omega.$$
Moreover,
\begin{gather}\label{y1}
|u|_{C^{1+\alpha, \frac{1+\alpha}{2}}(\overline{Q})}\leq 
C(\widehat{L})\Big(|y_0|_{L^2(\Omega)}
+|\hat\rho y_{1, d}|_{L^2(Q)}+|\hat\rho y_{2, d}|_{L^2(Q)}\Big).
\end{gather}
\end{lemma}

For the system (\ref{w1}), in Proposition $\ref{prop4}$ we consider
\begin{eqnarray*}
	&&b^{i j}=a^{i j}(z, \nabla z),\ (f_1, \cdots, f_n)^\top=F_2(z, \nabla z), \ f_0=F_1(z, \nabla z),\\[2mm]
	&&B^{i j}=A^{i j}(z, \nabla z),\  g_j=e^j(z, \nabla z)\ \mbox{ and } \ g_0=d_0(z, \nabla z).
\end{eqnarray*}
By the definition of $\mathcal{L}$,   we take
$$
\widehat L=\sum\limits_{i, j=1}^n |a^{i j}|_{C^2(\overline{B}_{n})}+
|f|_{C^2(\overline{B}_{n})},
$$
for $\overline{B}_{n}\deq \big\{\ (z, \zeta)\in\mathbb R^{n+1}\ \big|\ |z|\leq 1\ \mbox{ and }\ 
|\zeta|\leq 1\ \big\}.$ 
Further,  by (\ref{mu}), choose 
\begin{equation}\label{mu1}
	\hat\mu_0=Ce^{C\lambda^2}+\mu^*,
\end{equation}
for $\lambda=C(\widehat{L})e^{2\mu|\eta|_{C(\overline{\Omega})}}
=C(\widehat{L})e^{C\sum\limits_{i, j=1}^n (1+|a^{i j}|_{C^2(\overline{B}_{n})})|\eta|_{C(\overline{\Omega})}}$.

\medskip

\medskip

Now, we are in a position to give a proof of Theorem \ref{th1}.

\medskip

\noindent {\bf Proof of Theorem \ref{th1}. }
By Proposition \ref{DX!},  it suffices to prove the controllability result for the coupled  quasi-linear 
parabolic system (\ref{***}) (or (\ref{dl2})).

First, for any $z\in K$,  put
\begin{eqnarray*}
&&\Phi(z)=\Big\{\ y\in K\ \Big|\ \exists\ u\in 
C^{1+\alpha, \frac{1+\alpha}{2}}(\overline{Q}) \mbox{ and a positive constant } \widehat L,  \mbox{ such that } \\
&&\qquad\qquad\qquad\qquad (u, \bar y, p_1, p_2) \mbox{ satisfies } (\ref{w1}), (\ref{y1}),
 \mbox{ and } \bar{y}(x, T)=0 \mbox{ in } \Omega\   \Big\}.
\end{eqnarray*}
This defines a multi-valued  map $\Phi:K\rightarrow 2^K$,  provided that
$$|y_0|_{C^{3+\alpha}(\overline{\Omega})}
+|\hat\rho y_{1, d}|_{L^2(Q)}+|\hat\rho y_{2, d}|_{L^2(Q)}
+|y_{1, d}|_{C^{1+\alpha, \frac{1+\alpha}{2}}(\overline{Q})}
+|y_{2, d}|_{C^{1+\alpha, \frac{1+\alpha}{2}}(\overline{Q})} \ \mbox{  is small enough}.
$$
Indeed,  we apply the estimate in Lemma \ref{lemma*} for three equations in (\ref{w1}), respectively.
Since $|z|_{C^{3+\alpha, \frac{3+\alpha}{2}}(\overline{Q})}\leq 1$,   all coefficients of (\ref{w1})  
are in $C^{1+\alpha, \frac{1+\alpha}{2}}(\overline{Q})$, and \begin{eqnarray*}
&&\sum_{i, j, k=1}^n\Big([|a^{i j}(z, \nabla z)]_{x_k}|_{C^{1+\alpha, \frac{1+\alpha}{2}}(\overline{Q})}+
|[A^{i j}(z, \nabla z)]_{x_k}|_{C^{1+\alpha, \frac{1+\alpha}{2}}(\overline{Q})}\Big)\\
&&\quad+\sum_{i, j=1}^n\Big(|a^{i j}(z, \nabla z)|_{C^{1+\alpha, \frac{1+\alpha}{2}}(\overline{Q})}+
|A^{i j}(z, \nabla z)|_{C^{1+\alpha, \frac{1+\alpha}{2}}(\overline{Q})}\Big)+|d_0(z, \nabla z)|_{C^{1+\alpha, \frac{1+\alpha}{2}}(\overline{Q})}\\
&&\quad+\sum_{j=1}^n\Big(|e^j(z, \nabla z)|_{C^{1+\alpha, \frac{1+\alpha}{2}}(\overline{Q})}
+|F_2^j(z, \nabla z)|_{C^{1+\alpha, \frac{1+\alpha}{2}}(\overline{Q})}\Big)+|F_1(z, \nabla z)|_{C^{1+\alpha, \frac{1+\alpha}{2}}(\overline{Q})}\\
&&\leq C\big(n, T, \Omega, \alpha, \rho_0, |a^{i j}|_{C^4(\overline{B}_{n})}, 
|f|_{C^4(\overline{B_{n}})}\big).
\end{eqnarray*}
Hence, we obtain that
\begin{eqnarray*}
&&|\bar y|_{C^{3+\alpha, \frac{3+\alpha}{2}}(\overline{Q})}
\leq \widetilde C\Big(
|y_0|_{C^{3+\alpha}(\overline{\Omega})}
+|u|_{C^{1+\alpha, \frac{1+\alpha}{2}}(\overline{Q})}
+\frac{1}{\mu_1}|p_1|_{C^{1+\alpha, \frac{1+\alpha}{2}}(\overline{Q})}
+\frac{1}{\mu_2}|p_2|_{C^{1+\alpha, \frac{1+\alpha}{2}}(\overline{Q})}
\Big),\\[2mm]
&&|p_1|_{C^{3+\alpha, \frac{3+\alpha}{2}}(\overline{Q})}
\leq \widetilde C|\bar y-y_{1, d}|_
{C^{1+\alpha, \frac{1+\alpha}{2}}(\overline{Q})},\\[2mm]
&&\mbox{and }\ \  
|p_2|_{C^{3+\alpha, \frac{3+\alpha}{2}}(\overline{Q})}
\leq \widetilde C|\bar y-y_{2, d}|_
{C^{1+\alpha, \frac{1+\alpha}{2}}(\overline{Q})},
\end{eqnarray*}
where $\widetilde C$ depends on $n, T, \Omega, \alpha, \nu_1, \nu_2, \rho_0,  |a^{i j}|_{C^4(\overline{B}_{n})}$ and 
$|f|_{C^4(\overline{B_{n}})}.$
Adding  two inequalities on $p_1$ and $p_2$ together, we have 
\begin{eqnarray*}
&&|p_1|_{C^{3+\alpha, \frac{3+\alpha}{2}}(\overline{Q})}
+|p_2|_{C^{3+\alpha, \frac{3+\alpha}{2}}(\overline{Q})}\\[2mm]
&&\leq \widetilde C\Big(
|\bar y|_{C^{1+\alpha, \frac{1+\alpha}{2}}(\overline{Q})}
+|y_{1, d}|_{C^{1+\alpha, \frac{1+\alpha}{2}}(\overline{Q})}
+|y_{2, d}|_{C^{1+\alpha, \frac{1+\alpha}{2}}(\overline{Q})}\Big)\\[2mm]
&&\leq \widetilde C\Big(
|y_0|_{C^{3+\alpha}(\overline{\Omega})}+|u|_{C^{1+\alpha, \frac{1+\alpha}{2}}(\overline{Q})}
+|y_{1, d}|_{C^{1+\alpha, \frac{1+\alpha}{2}}(\overline{Q})}
+|y_{2, d}|_{C^{1+\alpha, \frac{1+\alpha}{2}}(\overline{Q})}\\[2mm]
&&\quad\quad\quad+\frac{1}{\mu_1}|p_1|
_{C^{1+\alpha, \frac{1+\alpha}{2}}(\overline{Q})}
+\frac{1}{\mu_2}|p_2|
_{C^{1+\alpha, \frac{1+\alpha}{2}}(\overline{Q})}\Big).
\end{eqnarray*}
When $\mu_k\geq 4\widetilde C$ for $k=1, 2$, it holds that
\begin{eqnarray*}
&&|p_1|_{C^{3+\alpha, \frac{3+\alpha}{2}}(\overline{Q})}
+|p_2|_{C^{3+\alpha, \frac{3+\alpha}{2}}(\overline{Q}))}\\[2mm]
&&\leq \widetilde C\Big(
|y_0|_{C^{3+\alpha}(\overline{\Omega})}+|u|_{C^{1+\alpha, \frac{1+\alpha}{2}}(\overline{Q})}
+|y_{1, d}|_{C^{1+\alpha, \frac{1+\alpha}{2}}(\overline{Q})}
+|y_{2, d}|_{C^{1+\alpha, \frac{1+\alpha}{2}}(\overline{Q})}\Big).
\end{eqnarray*}
Therefore, by (\ref{y1}), it follows that
\begin{eqnarray*}
&&|\bar y|_{C^{3+\alpha, \frac{3+\alpha}{2}}(\overline{Q})}+
|p_1|_{C^{3+\alpha, \frac{3+\alpha}{2}}(\overline{Q})}
+|p_2|_{C^{3+\alpha, \frac{3+\alpha}{2}}(\overline{Q})}\\[2mm]
&&\leq  \widetilde C C(\widehat{L})\Big(
|y_0|_{C^{3+\alpha}(\overline{\Omega})}
+|y_{1, d}|_{C^{1+\alpha, \frac{1+\alpha}{2}}(\overline{Q})}
+|y_{2, d}|_{C^{1+\alpha, \frac{1+\alpha}{2}}(\overline{Q})}
+|\hat\rho y_{1, d}|_{L^2(Q)}+|\hat\rho y_{2, d}|_{L^2(Q)}
\Big).
\end{eqnarray*}
This shows that $\Phi:K\rightarrow 2^K$ for a sufficiently small $\delta_0$, if
$$
|y_0|_{C^{3+\alpha}(\overline{\Omega})}
+\sum_{k=1}^2|y_{k, d}|_{C^{1+\alpha, \frac{1+\alpha}{2}}(\overline{Q})}
+\sum_{k=1}^2|\hat\rho y_{k, d}|_{L^2(Q)}\leq \delta_0.
$$ 
 Moreover,  we choose 
\begin{equation}\label{mu**}
\mu_0=4\widetilde C+\hat\mu_0 \quad(\mbox{recall } \hat\mu_0 \mbox{ in } (\ref{mu1})).
\end{equation}

Similar to the same arguments of Theorem 1.1 in  \cite{Liu} and  \cite{Liux}, by  Kakutani's fixed point theorem, there exists a $\bar y \in K$ such that $\bar y \in \Phi(\bar y)$.
This means that for the system $(\ref{dl2})$, there exists a control $\bar{u}\in 
C^{1+\alpha, \frac{1+\alpha}{2}}(\overline{Q})$, 
such that the corresponding solution to $(\ref{dl2})$ satisfies $\bar y(x, T)=0$ in $\Omega$.
\endpf

\subsection{Proof of Theorem \ref{+6}}\label{subsection}

In this subsection,  we  give a proof of Theorem \ref{+6}. By Theorem \ref{th1}, 
it suffices to prove that under suitable conditions,  for any leader control $u$, a Nash quasi-equilibrium for $J_1$ and $J_2$ is indeed a 
Nash equilibrium.

\medskip

\noindent {\bf Proof of Theorem \ref{+6}.} For any leader control $u\in B_{\rho_1}$, any $v_1\in B_{\rho_1}$, and
 the associated Nash quasi-equilibrium $(\bar v_1, \bar v_2)\in  (B_{\rho_1})^2$ for the 
 functionals $J_1$ and $J_2$, 
  by the Taylor expansion,  it holds that
  \begin{eqnarray*}
  &&J_1(v_1, \bar v_2; u)-J_1(\bar v_1, \bar v_2; u)\\[2mm]
  &&=J_{1, v_1}(\bar v_1, \bar v_2; u)(v_1-\bar v_1)
  +\int^1_0 (1-s) J_{1, v_1 v_1}\big((1-s)\bar v_1+sv_1, \bar v_2; u\big)(v_1-\bar v_1)^2ds,
  \end{eqnarray*}
where  $J_{1, v_1 v_1}$ denotes second-order  G\^ateaux derivative operator of $J_1$ 
with respect to the first variable, and for 
$v\in B_{\rho_1}$ and $w=v_1-\bar v_1$,
$$
\displaystyle J_{1, v_1 v_1}(v, \bar v_2; u)w^2
\deq\lim\limits_{h\rightarrow 0}\frac{1}{h}
\Big[J_{1, v_1}(v+hw, \bar v_2; u)w-J_{1, v_1}(v, \bar v_2; u)w\Big].
$$
Since $(\bar v_1, \bar v_2)$ is the Nash quasi-equilibrium, it follows that 
$J_{1, v_1}(\bar v_1, \bar v_2; u)w=0$. 
Hence, in order to prove  $(\bar v_1, \bar v_2)$ to be the Nash equilibrium for $J_1$ and $J_2$, 
it suffices to prove that 
\begin{equation}\label{DD1}
J_{1, v_1 v_1}(v, \bar v_2; u)w^2\geq 0\ \mbox{ and }\ 
J_{2, v_2 v_2}(\bar v_1, v; u)w^2\geq 0,\quad \forall\ v\in B_{\rho_1},
\end{equation}
where $$
\displaystyle J_{2, v_2 v_2}(\bar v_1, v; u)w^2
\deq\lim\limits_{h\rightarrow 0}\frac{1}{h}
\Big[J_{2, v_2}(\bar v_1, v+hw; u)w-J_{2, v_2}(\bar v_1, v; u)w\Big].
$$

In the following, the proof of (\ref{DD1}) is divided into two parts.

\medskip

{\bf Step 1. } For any $v\in B_{\rho_1}$,  we have that
\begin{eqnarray}\label{DD2}
\begin{array}{rl}
&\displaystyle J_{1, v_1}(v, \bar v_2; u)w
=\lim\limits_{h\rightarrow 0}\frac{J_1(v+hw, \bar v_2; u)-J_1(v, \bar v_2; u)}{h}\\[2mm]
&\displaystyle=\mu_1 \int^T_0\int_{\omega_1} vw dxdt+
\nu_1\int_Q \xi_*(y-y_{1, d})p dxdt,
\end{array}
\end{eqnarray}
where $y$ and $p$ satisfy the following system:
\begin{eqnarray}\label{DD3}
\left\{
\begin{array}{ll}
y_{t}-\sum\limits^n_{i, j=1}\big(a^{ij}(y,\nabla  y)y_{ x_i}\big)_{x_j}
+f(y, \nabla y)=\xi_{0}u+
\xi_{1} v
+\xi_{2} \bar v_2     &\mbox{ in }   Q,\\
p_{t}-\sum\limits^n_{i, j=1} \big(A^{ij}(y, \nabla  y)
p_{x_i}\big)_{x_j}+\sum\limits_{j=1}^n 
\big(e^j(y, \nabla y)  p\big)_{x_j}-d_0(y, \nabla  y)
p=
\xi_1w &\mbox{ in }Q,\\
y=p=0  &\mbox{ on }  \Sigma,\\[2mm]
y(x, 0)=y_0(x), \ p(x,  0)=0 &\mbox{ in }  \Omega.
\end{array}
\right.
\end{eqnarray}
By the Schauder estimates and energy estimates for parabolic equations of second order, 
we have  that
\begin{eqnarray}\label{DD4}
\begin{array}{rl}
&\displaystyle |y|_{C^{2+\alpha, 1+\frac{\alpha}{2}}(\overline{Q})}
\leq C\Big(|u|_{C^{\alpha, \frac{\alpha}{2}}(\overline{Q})}
+|v|_{C^{\alpha, \frac{\alpha}{2}}(\overline{Q})}
+|\bar v_2|_{C^{\alpha, \frac{\alpha}{2}}(\overline{Q})}
+|y_0|_{C^{2+\alpha}(\overline{\Omega})}
\Big)\\[4mm]
&\quad\displaystyle\leq C\Big(|u|_{C^{\alpha, \frac{\alpha}{2}}(\overline{Q})}
+\rho_1
+|\bar v_2|_{C^{\alpha, \frac{\alpha}{2}}(\overline{Q})}
+|y_0|_{C^{2+\alpha}(\overline{\Omega})}
\Big),\\[4mm]
&\displaystyle |p|_{C^{2+\alpha, 1+\frac{\alpha}{2}}(\overline{Q})}
\leq C(n, T, \Omega, \alpha, a^{i j}, f, \rho_1, u, \bar v_2, y_0)|w|_{C^{\alpha, \frac{\alpha}{2}}(\overline{Q})},\\[4mm]
&\mbox{and }\ \ \displaystyle |p|_{L^2(Q)}+|\nabla p|_{L^2(Q)}\leq C(n, T, \Omega, \alpha, a^{i j}, f, 
\rho_1, u, \bar v_2, y_0)|\xi_1w|_{L^2(Q)}.
\end{array}
\end{eqnarray}

\medskip
 
Furthermore,  by (\ref{DD2}),  we have
\begin{eqnarray}\label{DD5}
\begin{array}{rl}
&\displaystyle J_{1, v_1 v_1}(v, \bar v_2; u)w^2
=\lim\limits_{h\rightarrow 0}\frac{1}{h}
\Big[J_{1, v_1}(v+hw, \bar v_2; u)w-J_{1, v_1}(v, \bar v_2; u)w\Big], \\[2mm]
&\displaystyle=\mu_1\int^T_0\int_{\omega_1} w^2 dxdt
+\lim\limits_{h\rightarrow 0}\frac{\nu_1}{h}
\Big[\int_Q\xi_* (y^h-y_{1, d})p^h dxdt-\int_Q\xi_*(y-y_{1, d})p dxdt\Big],
\end{array}
\end{eqnarray}
where $y^h$ and $p^h$ satisfy the following system:
\begin{eqnarray}\label{DD6}
\left\{
\begin{array}{ll}
y^h_{t}-\sum\limits^n_{i, j=1}\big(a^{ij}(y^h,\nabla  y^h)y^h_{ x_i}\big)_{x_j}
+f(y^h, \nabla y^h)=\xi_{0}u+
\xi_{1} (v+hw)
+\xi_{2} \bar v_2     &\mbox{ in }   Q,\\
p^h_{t}-\sum\limits^n_{i, j=1} \big(A^{ij}(y^h, \nabla  y^h)
p^h_{x_i}\big)_{x_j}+\sum\limits_{j=1}^n 
\big(e^j(y^h, \nabla y^h)  p^h\big)_{x_j}-d_0(y^h, \nabla  y^h)
p^h=
\xi_1w &\mbox{ in }Q,\\
y^h=p^h=0  &\mbox{ on }  \Sigma,\\[2mm]
y^h(x, 0)=y_0(x), \ p^h(x,  0)=0 &\mbox{ in }  \Omega.
\end{array}
\right.
\end{eqnarray}
By the Schauder estimates  for  parabolic equations of second order again, 
we have  that
\begin{eqnarray}\label{DD7}
\begin{array}{rl}
&\displaystyle |y^h|_{C^{2+\alpha, 1+\frac{\alpha}{2}}(\overline{Q})}\\[2mm]
&\displaystyle\leq C\Big(|u|_{C^{\alpha, \frac{\alpha}{2}}(\overline{Q})}
+|v|_{C^{\alpha, \frac{\alpha}{2}}(\overline{Q})}
+|w|_{C^{\alpha, \frac{\alpha}{2}}(\overline{Q})}
+|\bar v_2|_{C^{\alpha, \frac{\alpha}{2}}(\overline{Q})}
+|y_0|_{C^{2+\alpha}(\overline{\Omega})}
\Big)\\[4mm]
&\displaystyle\leq C\Big(|u|_{C^{\alpha, \frac{\alpha}{2}}(\overline{Q})}
+3\rho_1
+|\bar v_2|_{C^{\alpha, \frac{\alpha}{2}}(\overline{Q})}
+|y_0|_{C^{2+\alpha}(\overline{\Omega})}
\Big),\\[4mm]
&\mbox{and }\ \ \displaystyle |p^h|_{C^{2+\alpha, 1+\frac{\alpha}{2}}(\overline{Q})}
\leq C(n, T, \Omega, \alpha, a^{i j}, f, \rho_1, u, \bar v_2, y_0)|w|_{C^{\alpha, \frac{\alpha}{2}}(\overline{Q})}.
\end{array}
\end{eqnarray}

\medskip

{\bf Step 2. } Introduce the following auxiliary system:
\begin{eqnarray}\label{DD8}
\left\{
\begin{array}{ll}
-q^h_{t}-\sum\limits^n_{i, j=1}\big(A^{ij}(y^h,\nabla  y^h)q^h_{ x_i}\big)_{x_j}
-\sum\limits_{j=1}^n 
e^j(y^h, \nabla y^h)  q^h_{x_j}\\[2mm]
\quad\quad\quad\quad\quad\quad\quad\quad\quad\quad\quad\quad\quad\quad-d_0(y^h, \nabla  y^h)
q^h=\xi_*(y^h-y_{1, d})     &\mbox{ in }   Q,\\[2mm]
-q_{t}-\sum\limits^n_{i, j=1}\big(A^{ij}(y,\nabla  y)q_{ x_i}\big)_{x_j}
-\sum\limits_{j=1}^n 
e^j(y, \nabla y)  q_{x_j}\\[2mm]
\quad\quad\quad\quad\quad\quad\quad\quad\quad\quad\quad\quad\quad\quad
-d_0(y, \nabla  y)
q=\xi_*(y-y_{1, d})     &\mbox{ in }   Q,\\[2mm]
q^h=q=0  &\mbox{ on }  \Sigma,\\[2mm]
q^h(x, T)=q(x,  T)=0 &\mbox{ in }  \Omega,
\end{array}
\right.
\end{eqnarray}
where $y^h$ and $y$ are the solutions  to (\ref{DD6}) and (\ref{DD3}).
Then, by (\ref{DD4}), 
\begin{eqnarray*}
&&|q|_{C^{2+\alpha, 1+\frac{\alpha}{2}}(\overline{Q})}
\leq C(n, T, \Omega, \alpha, a^{i j}, f, \rho_1, u, \bar v_2, y_0)(|y|_{C^{\alpha, \frac{\alpha}{2}}(\overline{Q})}+
|y_{1, d}|_{C^{\alpha, \frac{\alpha}{2}}(\overline{Q})})\\[2mm]
&&\leq C(n, T, \Omega, \alpha, a^{i j}, f, \rho_1, y_{1, d}, u, \bar v_2, y_0).
\end{eqnarray*}

On the other hand, by the duality between $p^h$ and $q^h$, and between $p$ and $q$,  we obtain that
$$
\displaystyle\int_Q \xi_* (y^h-y_{1, d})p^h dxdt=\int_Q \xi_1 q^h wdxdt \ 
\mbox{ and } \ 
\displaystyle\int_Q \xi_* (y-y_{1, d})p dxdt=\int_Q \xi_1 q wdxdt.
$$

By the above equalities and (\ref{DD5}), it holds that
\begin{eqnarray}\label{DD9}
\begin{array}{rl}
&\displaystyle J_{1, v_1 v_1}(v, \bar v_2; u)w^2
=\displaystyle\mu_1\int^T_0\int_{\omega_1} w^2 dxdt
+\lim\limits_{h\rightarrow 0}\frac{\nu_1}{h}
\int_Q \xi_1 w(q^h-q) dxdt\\[2mm]
&\displaystyle=\displaystyle\mu_1\int^T_0\int_{\omega_1} w^2 dxdt
+
\nu_1\int_Q \xi_1 wWdxdt,
\end{array}
\end{eqnarray}
where $W$ satisfies the following equation: 
\begin{eqnarray}\label{DD10}
\left\{
\begin{array}{ll}
\displaystyle-W_{t}-\sum\limits^n_{i, j=1}\big(A^{ij}(y,\nabla  y)W_{ x_i}\big)_{x_j}
-\sum\limits_{j=1}^n 
e^j(y, \nabla y)  W_{x_j}-d_0(y, \nabla y)W&\\[2mm]
\displaystyle\quad 
-\sum\limits^n_{i, j=1}\big(A^{ij}_y(y,\nabla  y)pq_{ x_i}\big)_{x_j}
-\sum\limits^n_{i, j=1}\big(\nabla_\zeta A^{ij}(y,\nabla  y)\cdot \nabla p q_{ x_i}\big)_{x_j}&\\[2mm]
\displaystyle\quad\quad
-\sum\limits_{j=1}^n 
e^j_y(y, \nabla y)pq_{x_j}-\sum\limits_{j=1}^n 
\nabla_\zeta e^j(y, \nabla y)\cdot\nabla p q_{x_j}&\\[3mm]
\displaystyle\quad\quad\quad-d_{0, y}(y, \nabla y)pq-\nabla_\zeta d_0(y, \nabla y)\cdot\nabla p q=\xi_* p
 &\mbox{ in }   Q,\\[2mm]
W=0  &\mbox{ on }  \Sigma,\\[2mm]
W(x, T)=0 &\mbox{ in }  \Omega.
\end{array}
\right.
\end{eqnarray}
We multiply the first equation of (\ref{DD10}) by 
$W$ and integrate the equality in $Q$. Notice that
$$
|y|_{C^{2+\alpha, 1+\frac{\alpha}{2}}(\overline{Q})}
+|q|_{C^{2+\alpha, 1+\frac{\alpha}{2}}(\overline{Q})}\leq C(n, T, \Omega, \alpha, a^{i j}, f, \rho_1, y_{1, d},
 u, \bar v_2, y_0)
$$
and 
$$
|p|_{L^2(Q)}+|\nabla p|_{L^2(Q)}\leq C(n, T, \Omega, \alpha, a^{i j}, f, \rho_1, u, \bar v_2, y_0)|w|_{L^2(\omega_1\times(0, T))}.
$$
Then, by the energy estimates,  it is easy to show that $$|W|_{L^2(Q)}\leq C(n, T, \Omega, \alpha, a^{i j}, f, \rho_1, y_{1, d}, u, \bar v_2, y_0)|w|_{L^2(\omega_1\times(0, T))}
\deq C_2|w|_{L^2(\omega_1\times(0, T))}.$$
By (\ref{DD9}),
\begin{eqnarray*}\label{DD*}
\begin{array}{rl}
&\displaystyle J_{1, v_1 v_1}(v, \bar v_2; u)w^2
=\displaystyle\mu_1\int^T_0\int_{\omega_1} w^2 dxdt
+
\nu_1\int_Q \xi_1 wWdxdt\\[2mm]
&\displaystyle\geq (\mu_1- C_2\nu_1)\int^T_0\int_{\omega_1} w^2 dxdt.
\end{array}
\end{eqnarray*}
Hence, when  we choose 
\begin{equation}\label{mu****}
\mu^0=\mu_0+C_2(\nu_1+\nu_2). \quad(\mbox{recall }  \mu_0 \mbox{ in } (\ref{mu**})), 
\end{equation}
 $J_{1, v_1 v_1}(v, \bar v_2; u)w^2\geq 0$, for any 
$v\in B_{\rho_1}$.

Similarly, we can prove that 
$J_{2, v_2 v_2}(\bar v_1,  v_2; u)w^2\geq 0$, for any 
$v\in B_{\rho_1}$.  This finishes the proof of Theorem \ref{+6}. \endpf

\end{document}